\documentclass [12pt] {report}
\usepackage[dvips]{graphicx}
\usepackage{epsfig, latexsym,graphicx}
\usepackage{amssymb}
\newcommand {\D}[2] {\displaystyle\frac{\partial{#1}}{\partial{#2}}}

\newcommand {\Dd}[3] {\displaystyle\frac{\partial^2{#1}}{\partial{#2}\partial{#3}}}

\newcommand {\la} {\lambda}

\newcommand {\de} {\delta}

\newcommand {\fr} {\displaystyle\frac}

\newcommand {\be} {\begin{equation}}
\newcommand {\ee} {\end{equation}}
\newcommand {\ba} {\begin{array}}
\newcommand {\ea} {\end{array}}
\newcommand {\bp} {\begin{picture}}
\newcommand {\ep} {\end{picture}}
\newcommand {\bc} {\begin{center}}
\newcommand {\ec} {\end{center}}
\newcommand {\bt} {\begin{tabular}}
\newcommand {\et} {\end{tabular}}
\newcommand {\lf} {\left}
\newcommand {\rg} {\right}

\newcommand {\cC} {{\cal C}}

\newcommand {\cI} {{\cal I}}

\newcommand {\cR} {{\cal R}}
\newcommand {\cS} {{\cal S}}

\newcommand {\ses} {\medskip}

\newcommand {\bibit} {\bibitem}
\newcommand {\nin} {\noindent}

\newcommand {\sign} {\mathop{\rm sign}\nolimits}

\newcommand {\wh} {\widehat}

 \usepackage{amsmath}

\def\2#1#2#3{{#1}_{#2}\hspace{0pt}^{#3}}
\def\3#1#2#3#4{{#1}_{#2}\hspace{0pt}^{#3}\hspace{0pt}_{#4}}
\newcounter{sctn}
\def\sec#1.#2\par{\setcounter{sctn}{#1}\setcounter{equation}{0}
                  \noindent{\bf\boldmath#1.#2}\bigskip\par}

    \newcommand {\s} {\ses}

\newcommand {\ff} {\vspace{1cm}}

\newcommand {\ed} {\end{document}}

\newcommand{\bvV}{\breve V}   \newcommand{\bvU}{\breve U}  \newcommand{\bvZ}{\breve Z}
\newcommand{\bvf}{\breve f}
\newcommand{\bvr}{\breve r}    \newcommand{\bvt}{\breve t}

\newcommand{\hL}{\hat L}

\newcommand{\hC}{\hat C}

\newcommand {\arcsinh} {\mathop{\rm arcsinh}\nolimits}

\newcommand {\FR}{F^{{\cal REG}}}

\newcommand {\FRR}{F^{{\cal REG-PSEUDO-FRD}}}

\newcommand {\FRRIN}  {F^{\cal REG-FRD-IN-PSEUDO-FRD}}

\begin {document}

\begin{titlepage}

\vspace{0.1in}

\begin{center}

{\bf \Large Two-axes pseudo-Finsleroid metrics: general  overview and angle-regular solution}

\end{center}

\vspace{0.3in}

\bc

\vspace{.15in}{\large G.S. Asanov\\}  \vspace{.25in}
{\it Devision of Theoretical Physics, Moscow State University\\
119992 Moscow, Russia\\
{\rm (}e-mail: asanov@newmail.ru{\rm)}} \vspace{.15in}

\ec

\begin{abstract}

\s  \s

The class of the  two-axes pseudo-Finslerian metrics which is specified by the condition of  the angle-separation in the  involved characteristic functions is proposed and  studied.
The complete    Total Set of algebraic and  differential equations
is derived in all rigor which are necessary and sufficient in order that a pseudo-Finsleroid metric function
belong to the class. It proves possible to solve the equations of the set.
The angle-regular solution  of   the Finsleroid-in-pseudo-Finsleroid type is found and described in detail.

\s

\nin

\s

\nin

\end{abstract}

\end{titlepage}

\vskip 1cm


\setcounter{equation}{0}

\vspace{0.3in}

\s

\nin
{\bf \large 1 Introduction}

\s \s

\nin
The applied ability of the methods of the Finsler Geometry is proportional to the variety
of  classes of the Finsler metric functions elaborated to reflect  various violations
of the spherical symmetry presupposed in the Euclidean and Riemannian Geometries.
Among possible violations,
the substitution of the
 axial symmetry with the spherical symmetry is the simplest case.
 Despite the  axial symmetry is widely appeared in theoretical and applied sciences,
in many patterns  the  axial symmetry is fulfilled but approximately.
The two-axes asymmetry    can well be regarded as the  much capable idea to proceed with.
Do the pseudo-Finsleroid metric functions possessing such an asymmetry exist when high regularity
properties are assumed?

We shall make a systematic  search to answer the question, following the method of introducing the angle dependence
and assuming the separation of angle dependencies, that was proposed and applied in the previous work [1-2].

Historically, Minkowski [3,4]  is well-known to have remarked that the Special Relativity
implies introducing the pseudo-Euclidean metric to geometrize the space-time.
After that, the pseudo-Euclidean geometry became raising its influence on the
philosophy as well as  computational methods of relativistic theories [5].
Beyond the square-root metric one may apply the ingenious methods of the Finsler geometry, including  the classical methods [6-9] as well as the recent and modern methods (see [9-15] and numerous references therein).


\s

{

We call a direction in the tangent space of a Finsler space {\it geometrically distinguished}
if the direction makes trace in the structure of the Finsler metric function.
Our goal in the present paper is to derive and develop the class of  four-dimensional pseudo-Finsleroid metric functions
which involve
two such directions, respectively of the vertical and horizontal geometrical meaning,
specifying the class by
 the   condition of angle-separated dependence
 of the characteristic functions.

{

Like to the preceding work  [1,2], we shall construct the pseudo-Finsleroid space on the four-dimensional
pseudo-Riemannian space, to be denoted by $\cR_4$.
To this end  we expand  the pseudo-Riemannian metric tensor $\{a_{ij}(x)\}$ of the space $\cR_4$
 with respect to an orthonormal vector frame
$\{b_i(x), i_i(x),j_i(x), i_{\{3\}i}(x)\}$, such that
\be
a_{ij}=b_ib_j-i_ii_j-j_ij_j-i_{\{3\}i}i_{\{3\}j}.
\ee
Respective 1-forms will be denoted by
$b=b_iy^i,i=i_iy^i,j=j_iy^i, i_{\{3\}}=i_{\{3\}i}y^i$.
The tensorial indices will be raised by means of the tensor  $\{a^{ij}(x)\}$
reciprocal to  $\{a_{ij}(x)\}$ , such that
 $a_{in}a^{nj}=\de^j_i, b^i=a^{in}b_n, {\it etc.}$
The consideration will be of local nature.
By $F=F(x,y)$ we shall denote a pseudo-Finslerian metric function.

{

\s \s

\nin
 {\bf  Definition 1.1} The {\it  pseudo-Finslerian cone} is
$$
\cC_x = \{  y\in T_xM :     F(x,y)=0\},
$$
which defines the hypersurface $\cC_x\subset T_xM$.
 If a vector  $y\in T_xM$  belongs to the cone, the vector is called {\it isotropic}.
  The vectors $y\in T_xM$ inside the  cone $\cC_x$ are called  $b$-{\it like}.

  \s  \s

\nin
This $\cC_x$  generalizes
 the conventional pseudo-Euclidean cone.

Below we confine our consideration to the $b$-like region of the tangent space, always assuming $b>0$.
Inside the cone $\cC_x$,  the function $F$ is assumed to be positive and  positively homogeneous with respect to tangent vectors, namely for any admissible $y\in T_xM$ and for an arbitrary   nonnegative
$t$ we have $F(x,ty)=tF(x,y)$.

\s  \s

\nin
 {\bf  Definition 1.2}
The hypersurface $\cI_x\subset T_xM$ introduced by
 $$
\cI_x = \{  y\in T_xM:  F(x,y)=1\}
$$
is called the {\it pseudo-Finslerian indicatrix}.
\ses

\s  \s

\nin
In the pseudo-Riemannian limit, the indicatrix $\cI_x$  reduces to the  pseudo-Euclidean
hyperboloid.

{

The straight line to which the vector $b^i(x)$ belongs in the tangent space $T_x$ will play the role of the {\it axis of the indicatrix}  $\cI_{x}$.
This axis  will be interpreted as {\it  vertical}.
Additionally, $i_{\{3\}}^i(x)$ will be interpreted as the vector which assigns
the {\it horizontal axis} of the indicatrix  $\cI_{x}$.

{

In Sect. 2, we
introduce the angle representations for the  Finslerian
 metric tensor $g_{ij}$ and  unit vector $l^i$.
We also construct the   orthonormal frame $\{l_i,u_i,m_i,p_i\}$,
 obtain   the tensor $C_{ijn}=(1/2)\partial g_{ij}/\partial y^n$,
and
propose the
  General Expansions of angle derivatives.
These objects  will play a fundamental role in the consideration performed in the subsequent sections.


In Sect. 3, after introducing
 the notion of the  separation of  angle dependence
 for the characteristic functions which enter the pseudo-Finsleroid metric function,
 we systematically derive  the Total Set of algebraic and differential equations
required to find solutions for the metric functions $F$ under study.
To this end,
the Particular Expansions of angle derivatives are derived which involve many simplifications
as compared to  the General  Expansions introduced in Sect. 2.

For existence of the functions $\eta=\eta(x,y), \theta=\theta(x,y),\phi=\phi(x,y)$ which enter
the Particular Expansions
the integrability  conditions   should be fulfilled which
involve   the    skew-vanishing lists $\cS_1, \cS_2,\cS_3$
which we shall  find.

The structural conditions of the studied class of the pseudo-Finsleroid metric functions
give rise to  the additional equations
written out in the First Group,  in the Second Group, and in the Third Group.

After that,
we   evaluate  the expansion of the tensor $C_{ijn}$
 in terms of the orthonormal vector frame. For the symmetry of the tensor,
the Symmetrizing Conditions should be set forth to specify the entered coefficients.
Formula (3.36) represents the resultant expansion,
which is truly convenient to use when evaluating the curvature tensor
$\wh R_{jpqn}=C^h{}_{pq}C_{hjn}-C^h{}_{pn}C_{hjq}$ in the tangent space.
Subjecting the tensor $\wh R_{jpqn}$ to  the requirement of the curvature constancy of indicatrix
entails the additional conditions, namely (3.39)-(3.41), and the value $\{-H^2\}$
for the curvature is
assigned by Formulas (3.43) and (3.44).
The inequality
$
H\ge 1
$
is always implied.
The Imperative Theorem is the completive assertion in  Sect. 3.
The angle separation  (3.1)-(3.5) introduced
complies with the separation presupposed  in [1,2].

{

In Sect. 4,
the   method of solving the equations entered the Total Set
is proposed and developed.
The possibility  to successfully proceed is to note that the equation triple (4.1)
entails the convenient ordinary first-order partial differential  equation (4.2) from which  the involved
key  function $\hL$
can readily be found (Formula (4.5)).
 The verification of  this implication is the one page evaluation, in the process of which
the differentiation of the second equality of (4.1) with respect to $\eta$ is the first step.
The first member of the triple  has been taken from the integrability condition  $\cS_1$, Formula  (3.14),
and two other members reflect the constancy of curvature of indicatrix (see Formulae (3.43)-(3.44)).
The obtained function $\hL$ involves one constant of integration,
which we  denote by $P$.

The solutions for $\hL$  are divided in three different classes: ~
{\it Class} I: $0<P<1$; ~
{\it Class} II: $P>1$, ~
{\it Class} III: $P<0$.

Additionally, the two constants $C$ and $C_7$ arise when the  separation of variables
is performed in the equations of
 the Second  Group and of  the   First  Group obtained in Sect. 3.

Thus the pseudo-Finsleroid space under consideration involves
three geometrically intrinsic scalars  $P=P(x),C=C(x), C_7=C_7(x)$
(which are  constants in each tangent space $T_x$).

The made observations  permit us to find dependence
of the characteristic functions $V$ and $r$ on $\eta$, according to Assertion 4.2 formulated to the end of Sect. 4.

{

In Sect. 5,
the  $\eta$-regular    solution
of the  $ {F^{\cal REG-FRD-IN-PSEUDO-FRD}}$-type will be proposed and described.
To this end
we choose  the {\it Class} II,
 assuming that $P>1$. We shall use the notation
$
\hC=1/{C},
$
where $C$ is the constant of
  Separation of Variables  that was arisen in the representation
$
\phi_t \phi_t=  C  Z_{tt}/Z
$
(see (4.9) in Sect. 4),
and introduce the constant $T$ by the help of the equality
$P=1/{T\hC}.
$
The appropriate idea is to subject the introduced constants to the inequalities
$
T>1, ~ ~ 0<\hC<1, ~ ~ T\hC<1.
$
The identification $C_7=1/P$ is made.


The dependence on  $\eta$ can  explicitly be found for the characteristic functions
$V=\breve V(x,\eta)$  and $r=\breve r(x,\eta)$. The obtained Formulas (5.4) and (5.14) representing the functions
$V=\breve V(x,\eta)$  and $r=\breve r(x,\eta)$
manifest  clearly the property
of $\eta$-regularity, namely that the functions
are smooth of class $C^{\infty}$ with respect to
$\eta$ over the total definition range $0 <\eta<\infty$.


This observation is geometrically of the {\it vertical meaning} and is not affected anyway by the geometry in the horizontal sections,
namely by the form of dependence of other characteristic functions
$\{U=\breve U(x,\theta), f=\breve f(x,\theta), Z=\breve Z(x,\phi),  t=\breve t(x,\phi)\}$
on $\theta$ and $\phi$.
However, it proves possible to subject the geometry to  the requirement  of constancy of curvature of the indicatrices in  the sectional horizontal spaces.
To this end, we assume the particular form (5.19) for the involved function $R_2$, which leads to success.
Indeed, the dependence of functions
$\{U=\breve U(x,\theta), f=\breve f(x,\theta)\}
$ on $\theta$
is found explicitly, according to (5.21) and (5.24).
Next,
 the simple  dependence   (see (5.28))
 of  the function $t=\breve t(x,\phi)$ on $\phi$ is proposed
 which fulfills the property
that the geometry is  $i_{\{3\}}$-axial, thereafter derivation  the  dependence of the function
$Z=\breve Z(x,\phi)$ on $\phi$ proves to be possible (see Formula (5.29)). With these solutions at hands, the attentive evaluation of the curvature of the
indicatrices in the sectional horizontal spaces lead to the striking  result that the curvature is constant and has the simple value (5.56).
The clear expression (5.57) is indicated for the determinant of the involved metric tensor of the horizontal sections,
which reveals the properties of  the  $\eta$-regularity and positive-definiteness.
The results of the evaluations performed are shortly summarized in Assertion 5.3,
which presents  the complete list (5.58) of the characteristic differential equations required.

{

In Conclusions, several important remarks have been made.

\s

 Appendix A has been added to show how  the Skew-Vanishing Lists of Sect. 3 can be verified
 by performing  attentive evaluation.

\s

Structurally, the evidence of the property of constant negative curvature for the indicatrix of the pseudo-Finsleroid space under study
is supported by the equalities (3.39)-(3.44) formulated in Sect. 3.
In Appendix B we derive the equalities upon evaluating the respective curvature tensor.


\s  \s

\setcounter{equation}{0}

\nin
{\bf \large 2  Angle representation}

\s  \s

\nin
We refer the consideration to the dimension
$
 N=4
$
and assume  the signature $\sign(g_{ij})=(+ - - -)$ for the Finslerian metric tensor $g_{ij}$.

\s  \s

\nin
{\bf  Primary Theorem} {\it Given a pseudo-Finsleroid metric function $F$.
The entailed Finslerian metric tensor admits the representation
\be
g_{ij}=
l_il_j
-
\fr{1}{H^2}
\lf(
\eta_i
\,\eta_j
+
\sinh^2\eta (
\theta_i\theta_j  +  \sin^2\theta  \phi_i\phi_j
)
\rg)
F^2,
\ee
where
$\eta(x,y),\theta(x,y),\phi(x,y)$ are three scalars homogeneous of  degree zero with respect to the argument $y,$
and $l_i=\partial F/\partial y^i.$
}

\s \s

\nin
The notation
 $\{\eta_i=\partial{\eta}/\partial{y^i}, \theta_i=\partial{\theta}/\partial{y^i}, \phi_i=\partial{\phi}/\partial{y^i}\}$, and
 $g_{ij}=(1/2)\partial F^2/\partial y^i\partial y^j$ is used.
The tensor
$
h_{ij}=g_{ij}- l_il_j
$
can be given by the formula
\be
h_{ij}=
-
\fr{1}{H^2}
\lf(
\eta_i
\,\eta_j
+
\sinh^2\eta (
\theta_i\theta_j  +  \sin^2\theta  \phi_i\phi_j
)
\rg)
F^2.
\ee


\s

The validity of the theorem can primarily
be recognized on the basis of  the differential geometry of indicatrix.
Such a geometry was systematically described in [7], applying the method of parametrical representation
of the indicatrix and using the respective projection factors. Namely, the Finslerian
indicatrix $I_x \subset T_x$ is considered as a hypersurface defined by the equation $F(x,y)=1$
in the tangent space $T_x$. In the present paper, we consider the $(N=4)$-dimensional case, so that the indicatrix is three-dimensional. Let $U_x\subset I_x$ be an open region on the indicatrix.
Fixing a point $x\in M$ and parameterizing the region by a set of three parameters $\{U^a\}=(U^1,U^2,U^3)$ assign the parametric representation
$l^i=t^i(x,U)$ to the unit tangent vectors $l^i=y^i/F$.
With  the help of the derivatives $t^i_a=\partial t^i/\partial U^a$
the indicatrix metric tensor $i_{ab}(x,U)$ is defined by the equality
$i_{ab}=-t^i_at^j_bh_{ij}$ (we have inserted the minus $\{-\}$ to respond to the signature (+ - - -) implied in the present paper). The parameters $U^a$ play the role of coordinates for the tensor
$i_{ab}$ which can geometrically be interpreted as a Riemannian metric tensor on the indicatrix.

Therefore, if we assume that the tensor $i_{ab}(x,U)$  corresponds to the constant negative
 curvature, we can use various known properties of the pseudo-Riemannian geometry. Among them, there exists
 the possibility of the choice   $\{U^1=\eta,U^2=\theta,U^3=\phi\}$ with  $\eta,\theta,\phi$ denoting
 hyperbolic angles.
With this choice, the components of the tensor $i_{ab}(x,U)$ are given by the list
\be
i_{11}=\fr{1}{H^2}, ~~  i_{22}=\fr{1}{H^2}\sinh^2\eta , ~~
i_{33}=\fr{1}{H^2}\sinh^2\eta \sin^2\theta ,
\ee
with the non-diagonal components vanished identically. The factor $1/H^2$ reflects the  property that
the  tensor corresponds to the curvature value $\{-H^2\}$, which we assume in the present paper.
Applying this method of parametrization of indicatrix,
the unit tangent vectors $l^i=y^i/F$ are given parametrically as
\be
l^i=t^i(x,\eta,\theta,\phi).
\ee


\s  \s

\nin
{\bf Definition 2.1} Formula (2.1) assigns the {\it angle representation of the Finslerian metric tensor} $g_{ij}$,
and  Formula (2.4)  introduces
the {\it angle representation of the unit vector} $l^i$.

\s  \s

Also, the square of the length element on the indicatrix
\be
ds^2=
  \fr1{H^2}
\lf((d\eta)^2+\sinh^2\eta\lf( (d\theta)^2+\sin^2\theta(d\phi)^2\rg) \rg)
\ee
is obtained.

From (2.1) we obtain the  expansion
\be
g_{ij}=
l_il_j- u_iu_j-m_im_j-p_ip_j,
\ee
where the {\it orthonormal frame} $\{l_i,u_i,m_i,p_i\}$  consists of the unit vector $l_i$ and the vectors
\be
u_i=\fr1HF\eta_i, ~~  m_i=\fr1HF\sinh\eta\, \theta_i,  ~~
p_i
=
\fr1H
F
\sinh\eta\,\sin\theta\,\phi_i.
\ee


By differentiation, from (2.1) we arrive at the following angle representation of tensor
$C_{ijn}=(1/2)\partial g_{ij}/\partial y^n$:

\s

$$
2C_{ijn}=
\fr1F(h_{in}l_j+l_ih_{jn})
-
\fr{2F}{H^2}
l_n
\Bigl(
\eta_i
\,\eta_j
+
                                                         \sinh^2\eta (
\theta_i\theta_j  +  \sin^2\theta  \phi_i\phi_j
)
\Bigr)
$$

\s \s

$$
-
\fr{F^2}{H^2}
\bigl(
\eta_{in}\eta_j  +  \eta_i\eta_{jn}
\bigr)
$$

\s \s

$$
-
\fr{F^2}{H^2}
\bigl(
\theta_{in}\theta_j  +  \theta_i\theta_{jn}
+
2\sin\theta  \cos\theta
\theta_n
\phi_i\phi_j
+
\sin^2\theta
(\phi_{in}\phi_j  +  \phi_i\phi_{jn})
\bigr)
\sinh^2\eta
$$

\s \s

$$
-
\fr{2F^2}{H^2}
\sinh\eta\cosh\eta
\eta_n (
\theta_i\theta_j  +  \sin^2\theta  \phi_i\phi_j
),
$$
or

\s \s

$$
2C_{ijn}=
\fr1F(h_{in}l_j+l_ih_{jn}-2l_nh_{ij})
-
\fr{F^2}{H^2}
\bigl(
\eta_{in}\eta_j  +  \eta_i\eta_{jn}
\bigr)
$$

\s \s

$$
-
\fr{F^2}{H^2}
\bigl(
\theta_{in}\theta_j  +  \theta_i\theta_{jn}
+
2\sin\theta  \cos\theta
\theta_n
\phi_i\phi_j
+
\sin^2\theta
(\phi_{in}\phi_j  +  \phi_i\phi_{jn})
\bigr)
\sinh^2\eta
$$

\s \s

\be
-
\fr{2F^2}{H^2}
\sinh\eta\cosh\eta
\eta_n (
\theta_i\theta_j  +  \sin^2\theta  \phi_i\phi_j
).
\ee


Since  the  orthonormal frame $\{l_i,u_i,m_i,p_i\}$ is complete
and the partial derivatives of $\eta_i,\theta_i,\phi_i$ with respect $y^j$ are symmetric, we can propose
the following {\it General Expansions of angle derivatives}:
$$
\eta_{ij}
=
-
\fr1F(l_i\eta_j+l_j\eta_i)   +   u_1\fr{1}{H^2}\sinh\eta \sin\theta(\phi_i\eta_j+\phi_j\eta_i)
+
u_2\fr{1}{H^2}
\sinh^2\eta \sin^2\theta \phi_i\phi_j  + u_3\fr{1}{H^2}\sinh^2\eta\theta_i\theta_j
$$

\s \s

\be
+
u_4\fr{1}{H^2}
\sinh\eta (\theta_i\eta_j+\theta_j\eta_i)
+
u_5\fr{1}{H^2}
\sinh^2\eta\sin\theta (\phi_i\theta_j+\phi_j\theta_i)
+
u_6\fr{1}{H^2}\eta_i\eta_j,
\ee

\s

\s

$$
\theta_{ij}=-\fr1F(l_i\theta_j+l_j\theta_i)
  +
  z_1\fr{1}{H^2}\sinh^2\eta\sin\theta (\phi_i\theta_j+\phi_j\theta_i)
 +
 z_2\fr{1}{H^2}
\sinh^2\eta\sin^2\theta  \phi_i\phi_j
+
z_3\fr{1}{H^2}\sinh^2\eta\theta_i\theta_j
$$

\s \s

\be
+
z_4\fr{1}{H^2}
\sinh\eta (\eta_i\theta_j+\eta_j\theta_i)
+
z_5\fr{1}{H^2}
\sinh\eta\sin\theta (\phi_i\eta_j+\phi_j\eta_i)
+
z_6\fr{1}{H^2}\eta_i\eta_j,
\ee
and

\s \s

$$
\phi_{ij}=-\fr1F(l_i\phi_j+l_j\phi_i)  +r_1\fr{1}{H^2}\sinh^2\eta\sin\theta(\phi_i\theta_j+\phi_j\theta_i)
 +
r_2\fr{1}{H^2}
\sinh^2\eta\sin^2\theta  \phi_i\phi_j
+
r_3\fr{1}{H^2}
\sinh^2\eta\theta_i\theta_j
$$

\s \s

\be
+
r_4\fr{1}{H^2}
\sinh\eta (\eta_i\theta_j+\eta_j\theta_i)
+
r_5\fr{1}{H^2}
\sinh\eta\sin\theta (\phi_i\eta_j+\phi_j\eta_i)
+
r_6\fr{1}{H^2}\eta_i\eta_j.
\ee
The coefficients $u_1,u_2,...,r_6$, which  are functions of $x$ and of the triple $\{\eta,\theta,\phi\}$
bear  the essential information on the Finsler metric functions that the expansions are referred to.
To be applicable, the expansions must be subjected to the conditions of integrability,
which require obviously that the skew-symmetrization of partial derivatives of the right-hand parts
in the expansions yields always zero.


\s  \s

\setcounter{equation}{0}

\nin
{\bf  \large 3 Separation of  angle dependence}

\s  \s

\nin
Henceforth, we assume for the  pseudo-Finslerian metric functions $F=F(x,y)$ of the considered type
 the  particular  structure specified by the condition of
{\it separation of  angle dependence}
\be
V=\bvV(x,\eta), ~ r=\bvr(x,\eta), ~ U=\bvU(x,\theta), ~ f=\bvf(x,\theta), ~ Z=\bvZ(x,\phi), ~ t=\bvt(x,\phi)
\ee
of the {\it characteristic functions}
 $\{V,r,U,f,Z,t\}$
which enter the  functions $F$ in accordance with the following
{\it Separation Scheme}:

\s \s

\be
F=b V(x,r),  ~~  r=z  U(x,f), ~~ f=c_2  Z(x,t)
\ee
with

\s

\be
 z=w_3, ~~ c_1=\fr {w_1}z, ~~  c_2=\fr {w_2}z,  ~~ ~~  w_1=\fr ib, ~~ w_2=\fr jb, ~~ w_3= \fr {i_{\{3\}}}b,
~~
 t=\fr{w_1}{w_2} =  \fr{c_1}{c_2}  =  \fr ij.
\ee
Here, $b$ is the 1-form introduced in Sect. 1.

  \s

It follows that
\be
\bvV(x,\eta)=V(x,\bvr(x,\eta)), ~~  \bvU(x,\theta)=U(x,\bvf(x,\theta)), ~~
\bvZ(x,\phi)=Z(x,\bvt(x,\phi)),
\ee
and

\s

\be
\eta=\bar\eta(x,V),  ~~ \theta=\bar\theta(x,f), ~~
\phi=\bar\phi(x,t).
\ee


The {\it structural conditions} (3.2)-(3.3) specify the dependence of the characteristic functions on
 tangent vectors $\{y^i\}$, and therefore the form of derivatives of the functions with respect to
$y^i$. In particular, we obtain the expansions
\be
l_i=b_iV+bV_i, ~~ V_i=V_{\eta}\eta_i, ~~ \eta_i=\eta_VV_i, ~~
\eta_{ij}=\eta_{VV}V_iV_j +  \eta_VV_{ij},
\ee
and

\s
\be
V_{i}=V_rr_i,  ~~  V_{ij}=V_{rr}r_ir_j+V_rr_{ij},
\ee
{\it etc}, where
the subscripts in $V,\eta,r$ mean differentiations; $l_i=\partial F/\partial y^i$.


The following assertion is of the fundamental  importance for our subsequent analysis.

\s  \s

\nin
{\bf Assertion 3.1} {\it
Assuming  the  separated  angle dependence} (3.1) {\it reduces
the  General Expansions} (2.9)-(2.11)  {\it  to the following  Particular Expansions
of angle derivatives:
}

\be
\eta_{ij}
=
-
\fr1F(l_i\eta_j+l_j\eta_i)
+
u_2\fr{1}{H^2}
\sinh^2\eta \sin^2\theta \phi_i\phi_j  + u_3\fr{1}{H^2}\sinh^2\eta\theta_i\theta_j
+
u_6\fr{1}{H^2}\eta_i\eta_j,
\ee

\s

\be
\theta_{ij}=-\fr1F(l_i\theta_j+l_j\theta_i)
  +
 z_2\fr{1}{H^2}
\sinh^2\eta\sin^2\theta  \phi_i\phi_j
+
z_3\fr{1}{H^2}\sinh^2\eta\theta_i\theta_j
+
z_4\fr{1}{H^2}
\sinh\eta (\eta_i\theta_j+\eta_j\theta_i),
\ee
and

\s

$$
\phi_{ij}=-\fr1F(l_i\phi_j+l_j\phi_i)  +r_1\fr{1}{H^2}\sinh^2\eta\sin\theta(\phi_i\theta_j+\phi_j\theta_i)
 +
r_2\fr{1}{H^2}
\sinh^2\eta\sin^2\theta  \phi_i\phi_j
$$

\s \s

\be
+
r_5\fr{1}{H^2}
\sinh\eta\sin\theta (\phi_i\eta_j+\phi_j\eta_i).
\ee

\s

\s

\nin
The verification of the validity of the assertion, namely that
the  separation of  angle dependence
 implies the nullifications
\be
u_1=u_4=u_5=z_1=z_5=z_6=r_3=r_4=r_6= 0
\ee
 in (2.10)-(2.12),
 involves calculations which are lengthy, although simple.
They are not reproduced in the present paper.

{

The following assertion is valid.

\s  \s

\nin
{\bf Assertion 3.2} {\it For existence of the functions $\eta=\eta(x,y), \theta=\theta(x,y),\phi=\phi(x,y)$ which obey
the set} (3.8)-(3.10) {\it of    Particular Expansions it is necessary and sufficient that the functions
fulfill the conditions $\cS_1, \cS_2,\cS_3$ listed below.
}

\s  \s


The {\it Skew-Vanishing Lists}   $\cS_1, \cS_2,\cS_3$ are
obtainable respectively from the integrability conditions
\be
\D{\eta_{ij}}{y^n}-\D{\eta_{in}}{y^j}=0, ~~  \D{\theta_{ij}}{y^n}-\D{\theta_{in}}{y^j}=0, ~~  \D{\phi_{ij}}{y^n}-\D{\phi_{in}}{y^j}=0.
\ee
Namely, when  the general expansions
\s
$$
u_{2n}=u_{2\eta}\eta_n+u_{2\theta}\theta_n+u_{2\phi}\phi_n, ~~
u_{3n}=u_{3\eta}\eta_n+u_{3\theta}\theta_n+u_{3\phi}\phi_n, ~~
u_{6n}=u_{6\eta}\eta_n+u_{6\theta}\theta_n+u_{6\phi}\phi_n
$$
are inserted in (3.8),
we come to the  following list.

\s

\s

\nin
{\bf The  Skew-Vanishing List   $\cS_1$: }

\s

\be
u_{2\theta}
=
u_{2\phi}
=
u_{6\theta}
=
u_{6\phi}
=0,
\ee

\s

\be
u_{2\eta}+\fr{1}{H^2}u_2(u_2-u_6)=1,
\ee
and
\be
z_{4\theta}
=
z_{4\phi}
=0, \quad
(r_5\sin\theta)_{\theta}
=
r_{5\phi}
=0.
\ee

\s  \s


\nin
Next,
with  the derivative expansions

\s
$$
z_{2n}=z_{2\eta}\eta_n+z_{2\theta}\theta_n+z_{2\phi}\phi_n,  ~~
z_{3n}=z_{3\eta}\eta_n+z_{3\theta}\theta_n+z_{3\phi}\phi_n,  ~~
z_{4n}=z_{4\eta}\eta_n,
$$
from (3.9) we obtain  the  following equations.

\s \s

\nin
{\bf The  Skew-Vanishing List   $\cS_2$: }

\s

\be
z_{3\phi}=0, ~~
z_{2\eta} \sinh\eta + 2z_2 \cosh\eta =0, ~~
z_{3\eta}\sinh\eta
+
2z_3
 \cosh\eta
=0,
\ee

\s

\be
z_{2\theta}
-
z_4u_2
\fr{1}{H^2}
\sinh\eta
+
\fr{1}{H^2}
z_2
(
z_2
-
z_3
)
\sinh^2\eta
=
1,
\ee
and

\s

\be
z_{4\eta}\sinh\eta + z_4\cosh\eta
+
\fr{1}{H^2}
z_4
\lf(
u_6
\sinh\eta
-z_4
\sinh^2\eta
\rg)
=
-1.
\ee

\s  \s

\nin
From (3.16) we can make the significant conclusion
\be
(z_2\sinh^2\eta)_{\eta}=(z_3\sinh^2\eta)_{\eta}=0.
\ee


Finally, applying the expansions

$$
r_{1n}=r_{1\eta}\eta_n+r_{1\theta}\theta_n+r_{1\phi}\phi_n, ~~
r_{2n}=r_{2\eta}\eta_n+r_{2\theta}\theta_n+r_{2\phi}\phi_n,  ~~
r_{5n}=r_{5\eta}\eta_n+r_{5\theta}\theta_n
$$
to (3.10) leads  to  the  following result.

\s  \s

\nin
{\bf The  Skew-Vanishing List   $\cS_3$: }

\s

\be
(r_2\sinh^2\eta)_{\eta}=0, \qquad  (r_2\sin^2\theta)_{\theta}=0,  ~~
 r_{1\phi}=0,
\ee

\s

\be
 r_{1\theta}
 \sin\theta
 +
 r_1
\lf(
-
r_1\fr{1}{H^2}\sinh^2\eta\sin^2\theta
+
\fr1{H^2}z_3
\sinh\eta\sin\theta
+
 \cos\theta
\rg)
+
r_5u_2
\sinh\eta\sin\theta
=
-1,
\ee
and

\s

\be
r_{5\eta}
\sinh\eta\sin\theta
+r_5
\lf(
u_6\fr{1}{H^2}
\sinh\eta\sin\theta
-
r_5\fr{1}{H^2}
\sinh^2\eta\sin^2\theta
+
\cosh\eta\sin\theta
\rg)
=
-1,
\ee
together with
\be
  (r_2\sin^2\theta)_{\theta}=0,  ~~
 r_{1\phi}=0.
\ee

\s

\s

\nin
The validity of these lists  $\cS_1, \cS_2,\cS_3$ is verified in Appendix A.


\s

Next, we can draw additional  valuable information by using the structural expansions (3.6)-(3.7) and their concomitants.

Namely, attentive calculations lead to the conclusion that  the  initial equation (3.8)
for $\eta_{ij}$ is equivalent to the following list of vanishing.

\s \s

\nin
{\bf  First Group  of structural equations}:

\ff

\be
\lf(
\eta_{VV}  -u_6\fr{1}{H^2} \eta_V\eta_V
+\fr2V
\eta_V
\rg)
V_rV_r
+
  \eta_V
V_{rr}
=0,
\ee

\ff

\be
u_2\fr{1}{H^2}
\sinh^2\eta \sin^2\theta
 \phi_t \phi_t
-
\eta_V
V_r
U_f
Z_{tt}
w_2
=0,
\ee

\ff

\be
 u_3\fr{1}{H^2}\sinh^2\eta
\theta_f\theta_f
-
 \eta_V
 zU_{ff}
V_r
=0.
\ee


The same method can be applied to the expansion (3.9) obtained for $\theta_{ij}$.
We arrive at the additional  equations.

\s

\nin
{\bf  Second Group  of structural equations}:

\s \s

\be
\theta_{ff}
-
\lf(
z_3\fr{1}{H^2}\sinh^2\eta
\,
\theta_f
-
\fr2U
U_f
\rg)
\theta_f
=0,
\ee

\s

\be
z_2\fr{1}{H^2}
\sinh^2\eta \sin^2\theta
 \phi_t \phi_t
-
\theta_f
c_2Z_{tt}
=0,
\ee

\s

\be
\lf(
z_4\fr{1}{H^2}
\sinh\eta
\,
\eta_V
-
\fr1V
\rg)
z
V_r
+
\fr1U
=0.
\ee

\s  \s

Lastly, considering
 the expansion  (3.10) for $
\phi_{ij}$
leads to new series of equations.

 \s

\nin
{\bf Third Group  of structural equations}:

\s \s

\be
\phi_{tt}
=
-
\fr2Z
Z_t \phi_t
+
r_2\fr{1}{H^2}
\sinh^2\eta\sin^2\theta
\phi_t\phi_t,
\ee

\s

\be
r_1\fr{1}{H^2}\sinh^2\eta\sin\theta
\theta_f
-
\fr1U
U_f
+
\fr1f
=0,
\ee

\s

\be
\lf(
r_5\fr{1}{H^2}
\sinh\eta\sin\theta
\eta_V
-
\fr1V
\rg)
r
V_r
+
1
=0.
\ee

\s \s

\nin
The verification of the First, Second,  and Third Groups of equations
is the straightforward process and is not displayed in the present paper.


Now, we can evaluate the tensor $C_{ijn}$ using the representation (2.8) prepared to convenience in Sect. 2.
The tensor should be given by a symmetric representation.
It can readily be verified that the following conditions should be set forth.

\s \s

\nin
{\bf  Symmetrizing Conditions   for the tensor $C_{ijn}$}:

\be
z_4\sinh^2\eta=u_3\sinh\eta
-
2H^2
\cosh\eta, \quad
r_5\sinh^2\eta\sin\theta=u_2\sinh\eta
-
2H^2
\cosh\eta,
\ee

\s

\s

\be
z_2
\sinh^2\eta
\sin\theta
=
r_1
\sinh^2\eta \sin^2\theta
+
2H^2
 \cos\theta.
\ee

\s \s

\nin
If we compare the last equality with (3.19),
we just conclude  that
\be
(r_1\sinh^2\eta)_{\eta}=0.
\ee


With the Symmetrizing Conditions assumed, simple calculations lead to the following result.

\s \s

\nin
{\bf Assertion 3.3} {\it The tensor $C_{ijn}$ when written in terms of the
orthonormal  vectors has the structure

$$
FHC_{ijn}=
-  u_6u_iu_ju_n
-
L_3
 (u_im_j m_n+u_nm_jm_i  +  u_jm_im_n)
-
L_2
(u_ip_j p_n  + u_np_i p_j+u_jp_ip_n)
$$

\s

\s

\s

\be
-
\sinh\eta
 z_3m_im_jm_n
-
L
 (m_jp_ip_n+m_n p_ip_j+m_i p_jp_n)
-
 r_2
\sinh\eta
\sin\theta
p_ip_jp_n,
\ee
\ses\\
where

\be
L_2
=
   u_2  - H^2 \fr{\cosh\eta}{\sinh\eta}   , ~~
L_3
=
  u_3  - H^2 \fr{\cosh\eta}{\sinh\eta}   , ~~
L=
\lf(
r_1
\sinh^2\eta \sin\theta
+
H^2
\fr{ \cos\theta}{\sin\theta}
\rg)
\fr1{\sinh\eta}.
\ee
}

\s  \s

\nin
The  orthonormal frame $\{l_i,u_i,m_i,p_i\}$ (see (2.7) in Sect. 2) has been used.


Finally,  we evaluate  the curvature tensor
\be
\wh R_{jpqn}=C^h{}_{pq}C_{hjn}-C^h{}_{pn}C_{hjq}
\ee
 of the tangent space and set forth the requirement that the tensor manifests the property that
 the indicatrix is a space of constant negative curvature.
In the process of evaluation it is getting clear that the property implies the equalities
\be
L_2=L_3,
\ee
which  is equivalent to
\be
u_3=u_2
\ee
(see (3.37)),
and
\be
L_2
(u_6-2L_2)
=
L(z_3\sinh\eta-L).
\ee
The result reads
\be
F^2\wh R_{jpqn}
 =
T^*
(h_{pq}h_{jn}-h_{pn}h_{jq})\equiv
-T^*
(h_{pn}h_{jq}-h_{pq}h_{jn}),
\ee
where
\be
T^*-\fr1{H^2}
\lf( (L_2)^2
+
L_2
(u_6-2L_2)
\rg)
\ee
(see Appendix B).

To have the value $-H^2$ for the indicatrix curvature we should take
\be
T^*=1-H^2.
\ee

Since
$
u_3=u_2
$
(see (3.40)),
we can conclude that the equality
\be
z_4=r_5\sin\theta
\ee
holds owing to   the symmetrizing condition (3.34), from which it follows that
the last equation (3.32) in the Third Group is the same as the third equation (3.29) in the Second Group.


From  (3.41) and (3.43) it follows that
\s
\be
Lz_3\sinh\eta
=
-
H^2T^*
- (L_2)^2 +L^2.
\ee
The function $L$ was introduced in
 (3.37).
Using here the equality
$
z_2
\sinh^2\eta
\sin\theta
=
r_1
\sinh^2\eta \sin^2\theta
+
2H^2
 \cos\theta
$
which presents in the Symmetrizing Condition
(3.34),
we obtain the convenient representation
$$
L\sinh\eta=
z_2
\sinh^2\eta
-
H^2
\fr{ \cos\theta}{\sin\theta}.
$$
Thus,
 in terms of
the new notation
\be
\breve z_2= z_2\sinh^2\eta, ~~  \breve z_3= z_3\sinh^2\eta,
\ee
the equation (3.46) takes on the form of the following $\{\breve z_2,\breve z_3\}$-{\it quadratic equation}:

\be
\lf(\breve z_2
-
H^2
\fr{ \cos\theta}{\sin\theta}
\rg)
^2
-
\lf(\breve z_2
-
H^2
\fr{ \cos\theta}{\sin\theta}
\rg)  \breve z_3
-
\lf(H^2T^* + (L_2)^2\rg)\sinh^2\eta
=0.
\ee
Since
$
\breve z_{2\eta}=\breve z_{3\eta}=0
$
(see (3.19)), $\breve z_2$ and $\breve z_3$ are functions of  $\theta$
(and are independent of $\eta$ or $\phi$).

\s

The derivative of $\breve z_2$ can be obtained from (3.17), namely
$$
\breve z_{2\theta}
-
z_4\sinh^2\eta
u_2
\fr{1}{H^2}
\sinh\eta
+
\fr{1}{H^2}
\breve z_2
(
\breve z_2
-
\breve z_3
)
=
\sinh^2\eta.
$$
The Symmetrizing Conditions
(3.33) involves
$
z_4\sinh^2\eta=u_2\sinh\eta
-
2H^2
\cosh\eta.
$
We eventually obtain

\be
\breve z_{2\theta}
+
\fr{1}{H^2}
\breve z_2
(
\breve z_2
-
\breve z_3
)
=
\sinh^2\eta
+
\fr{1}{H^2}
\lf(u_2\sinh\eta
-
2H^2
\cosh\eta
\rg)
u_2
\sinh\eta.
\ee


To summarize up, we introduce the following definition.

\s

\nin
{\bf Definition 3.1}    The {\it Total Set} of the conditions and equations
involve the Skew-Vanishing integrability conditions $\cS_1, \cS_2,\cS_3$ and the Symmetrizing Conditions, supplemented by the structural equations
written out in the First Group,  in the Second Group, and in the Third Group,
and also  the following three   conditions (i), (ii), and (iii)  which  assure the constancy of indicatrix curvature:
(i)  $u_3=u_2$;  (ii)   $T^*=- ( (L_2)^2 + L_2 (u_6-2L_2) )/H^2 $ with $T^*=1-H^2$;
 (iii)  the $\{\breve z_2,\breve z_3\}$-quadratic equation must be fulfilled. The curvature value is equal to $-H^2$.

\s  \s

\nin
From the content of the present section we can make  the following general conclusion.

\s  \s

\nin
{\bf Imperative Theorem}   {\it In order that a   two-axes pseudo-Finsleroid
metric function $F$ be of the angle-separated  type, it is necessary and sufficient that the function be
compatible with each condition or equation which enters the  Total Set.
}

\s  \s


\setcounter{equation}{0}

\nin
{\bf  \large 4 Solving the derived equations}

\s \s

\nin
The  {\it Total Set} involves many  differential equations for numerous unknown functions. At the first sight,  the hope to resolve the set can be regarded as excessively optimistic.


However, a lucky and  simple possibility is offered when we pay due attention to the following equation triple:

\be
u_{2\eta}+\fr{1}{H^2}u_2(u_2-u_6)=1, ~~ L_2=   u_2  - H^2 \fr{\cosh\eta}{\sinh\eta}, ~~\fr1{H^2}  L_2  (L_2-u_6)=T^*\equiv 1-H^2.
\ee
Here, the first member has been taken from the integrability condition  $\cS_1$, Formula  (3.14),
and two other members reflect the constancy of curvature of indicatrix (see Formulae (3.40)-(3.44)).

\s

Indeed,  differentiation with respect to $\eta$ yields
$$
L_{2\eta}=
-\fr{1}{H^2}u_2(u_2-u_6)+1 +  H^2 \fr{1}{\sinh^2\eta}
=
-\fr{1}{H^2}u_2(L_2-u_6)
-u_2\fr{\cosh\eta}{\sinh\eta}
+1 +  H^2 \fr{1}{\sinh^2\eta},
$$
or

\s

$$
L_{2\eta}=
\fr{1}{H^2}L_2(u_6-L_2)
+ \fr{\cosh\eta}{\sinh\eta}
(u_6-L_2)
-u_2\fr{\cosh\eta}{\sinh\eta}
+1 +  H^2 \fr{1}{\sinh^2\eta},
$$
which can be written as
$$
L_{2\eta}=
\fr{1}{H^2}L_2(u_6-L_2)
+ \fr{\cosh\eta}{\sinh\eta}
(u_6-L_2)
-
\lf(  L_2+ H^2 \fr{\cosh\eta}{\sinh\eta}   \rg)
\fr{\cosh\eta}{\sinh\eta}
+1 +  H^2 \fr{1}{\sinh^2\eta}.
$$
We reduce this equality to read
$$
L_{2\eta}=
\fr{1}{H^2}L_2(u_6-L_2)
+ \fr{\cosh\eta}{\sinh\eta}
(u_6-L_2)
-
  L_2
\fr{\cosh\eta}{\sinh\eta}
+1 -  H^2.
$$
Multiplying by $L_2$ yields

\s

$$
L_2L_{2\eta}=
-L_2T^*
- \fr{\cosh\eta}{\sinh\eta}
H^2T^*
-
  (L_2)^2
\fr{\cosh\eta}{\sinh\eta}
-L_2(  H^2-1)
=
 \fr{\cosh\eta}{\sinh\eta}
H^2(H^2-1)
-
  (L_2)^2
\fr{\cosh\eta}{\sinh\eta}.
$$

\ff

In this way, we arrive at the validity of the following assertion.


\s  \s

\nin
{\bf Assertion 4.1} {\it The function $L_2$ should be obtained from the following simple differential equation:}
\be
L_2L_{2\eta}=\lf(H^2(H^2-1)-(L_2)^2\rg)\fr{\cosh\eta}{\sinh\eta}.
\ee

\s \s

{

\nin
The equation can readily be solved, yielding
\be
\lf( \fr1{H} L_2\rg)^2
=
H^2-1
-
\fr{H^2}{P\sinh^2\eta}
(1-P)
,
\ee
where $P$ is an integration constant,
or

\be
\lf(\fr1{H^2}  \sinh\eta L_2 \rg)^2
=\hL,
\ee
with
\be
\hL=1-\fr1{P}+\lf(1-\fr1{H^2}\rg)  \sinh^2\eta \equiv
\fr1{H^2}-\fr1{P}+\lf(1-\fr1{H^2}\rg)\cosh^2\eta.
\ee

\ff

The solutions are divided in three different classes: ~
{\it Class} I: $0<P<1$; ~
{\it Class} II: $P>1$, ~
{\it Class} III: $P<0$.


Let us evaluate the last term in  the
$\{\breve z_2,\breve z_3\}$-{quadratic equation} (3.48), making use of  (4.4). We get

$$
\fr1{H^4}  \sinh^2\eta
\lf[
H^2T^*+ (L_2)^2
\rg]
=
\fr1{H^2}T^*\sinh^2\eta +  \hL
=
$$

\s

$$
\fr1{H^2}(1-H^2)\sinh^2\eta   +  1-\fr1{P}+\lf(1-\fr1{H^2}\rg)  \sinh^2\eta
= 1-\fr1{P}
$$
So, the considered equation (3.48) can explicitly be written as
\be
\lf(\breve z_2
-
H^2
\fr{ \cos\theta}{\sin\theta}
\rg)
^2
-
\lf(\breve z_2
-
H^2
\fr{ \cos\theta}{\sin\theta}
\rg)  \breve z_3
-
\lf(1-\fr1P\rg)H^4
=0.
\ee
Using the notation
\be
\breve L_2=\breve z_2
-
H^2
\fr{ \cos\theta}{\sin\theta},
\ee
we can write the equation in the form
\be
\fr1{H^4}\breve L_2(\breve L_2-\breve z_3)
=
1-\fr1P,
\ee
which is quite similar to (3.43).

{

The {\it second line of simplifications} is opening up
when we refer to (3.28) of the Second  Group in Sect. 3 and write
$$
z_2\fr{1}{H^2}
\sinh^2\eta \sin^2\theta
\fr1{\theta_ff}
=
\fr1Z
\fr1{ \phi_t \phi_t}
Z_{tt},
$$
where we have used the equality  $c_2=f/Z$ which was written in (3.2) in Sec. 3.
The equation  is remarkable in that the left-hand part is a function of $\theta$ (it should be noted that
$(z_2\sinh^2\eta)_{\eta}=0$ according to (3.19) in Sect. 3), while the right-hand part
does not depend on $\theta$ and is a function of the argument  $t$.
Therefore, the equation implies the {\it  separation of variables}:
\be
\theta_ff
=
z_2\fr{1}{H^2}
C\sinh^2\eta \sin^2\theta, ~~  ~~
\phi_t \phi_t=
C
\fr {Z_{tt}}Z, \qquad C=C(x),
\ee
where $C$ is  the {\it first constant of  separation of variables}.

Introducing the function
\be
R_2=\fr{1}{H^2}
z_2\sinh^2\eta \sin\theta,
\ee
we get
\be
\theta_ff
=
R_2 C\sin\theta.
\ee

\s  \s


{

Keeping in mind the identification  $u_3=u_2$ (see (3.40) in Sect. 3), on compare (3.25) with (3.26) of the First  Group we can conclude

$$
\fr1{ \phi_t \phi_t}
U_f
Z_{tt}
w_2
=
\fr1{\theta_f\theta_f}
 zU_{ff} \sin^2\theta.
$$
Since $f=c_2Z$ and  $w_2=c_2z$, we obtain

$$
U_ff
\fr1{Z \phi_t \phi_t}
Z_{tt}
=
\fr1{\theta_f\theta_f}
U_{ff} \sin^2\theta,
$$
or

\s
\be
CU_{ff} \sin^2\theta=
U_ff
\theta_f\theta_f.
\ee


\s

Also, considering  (3.25) of the  First Group

$$
u_2\fr{1}{H^2}
\sinh^2\eta \sin^2\theta
=
\eta_V
V_r
U_f
\fr1{ \phi_t \phi_t}
Z_{tt}
w_2,
$$
we are coming to

\s

\be
u_2\fr{1}{H^2}
\sinh^2\eta \sin^2\theta
=
\eta_V
V_rr
\fr1U
U_ff
\fr1{C}.
\ee
From the previous equality
the {\it second angle-separation line} starts
\be
\fr{1}{H^2} u_2\sinh^2\eta
=
C_7\eta_rr, ~~
\fr1U
U_ff
\fr1{C}
=
C_7 \sin^2\theta, ~~~
C_7=C_7(x),
\ee
where the vanishing
 $u_{2\theta} = u_{2\phi} = 0 $
written in (3.12) has been taken into account.
$C_7$ is  the   {\it second constant of  separation of variables}.

{

With (4.14), the dependence of the function  $r=\breve r(x,\eta)$ on $\eta$ can be found upon integration:
\be
\ln r=C_7\int\fr{d\eta}{ \fr{1}{H^2} u_2
\sinh^2\eta
 }
\ee
with
$$
u_2=   L_2  + H^2 \fr{\cosh\eta}{\sinh\eta}
$$
and $L_2$ taken from (4.3).

{

To find also $V_r$, we
compare the third line (3.29) of the Second Group
with first line  (3.33) in  the Symmetrizing Condition.
We obtain the equation

\be
\fr{1}{H^2}
\lf(
u_2 - 2H^2 \fr{\cosh\eta}{\sinh\eta}
\rg)
\eta_rr
-
\fr1V
V_rr
+
1
=0,
\ee
where $r=zU$ has been used.
Multiplying by $C_7$ and using
$
(1/H^2) u_2 \sinh^2\eta = C_7\eta_rr $
(see (4.14)),
we conclude that

\s

$$
\fr{1}{H^2}
\lf(
u_2 - 2H^2 \fr{\cosh\eta}{\sinh\eta}
\rg)
\fr{1}{H^2} u_2
\sinh^2\eta
+
C_7(-
\fr1V
V_rr
+
1
)
=0,
$$
or

$$
\fr1V
V_rr
=
\fr{1}{C_7H^2}
\lf(
u_2 - 2H^2 \fr{\cosh\eta}{\sinh\eta}
\rg)
\fr{1}{H^2} u_2
\sinh^2\eta
+
1.
$$
Last, taking into account the equality

$$
u_2=   L_2  + H^2 \fr{\cosh\eta}{\sinh\eta},
$$
we  arrive  at the  representation

$$
\fr1V
V_rr
=
\fr{1}{C_7H^2}
\lf(
L_2 - H^2 \fr{\cosh\eta}{\sinh\eta}
\rg)
\lf(
L_2 + H^2 \fr{\cosh\eta}{\sinh\eta}
\rg)
\fr{1}{H^2}
\sinh^2\eta
+
1.
$$
Inserting here
$$
\fr1{H^2}
(L_2)^2
=
H^2-1
-
\fr{H^2}{P\sinh^2\eta}
(1-P)
$$
(see (4.3)) yields the following simple result:

\be
\fr1V
V_rr
=
-\fr{1}{C_7H^2}
\sinh^2\eta
-
\fr{1}{C_7P}+ 1.
\ee

{

It is convenient to use the function
\be
 R_1= \cosh\eta  +\sqrt{\hL},
 \ee
 such that
 \be
\fr1{H^2} u_2 \sinh\eta= R_1;
 \ee
$\hL$ is the function that was introduced in (4.4)-(4.5).
By the help of   (4.19) and (4.14),
we obtain the equations
\be
\eta_rr=
\fr1{C_7}
R_1\sinh\eta
\ee
and

\be
\fr1V
V_{\eta}
=
-
\fr1{R_1\sinh\eta}
\lf(
\fr{1}{H^2}
\sinh^2\eta
+
\fr{1}{P}- C_7\rg),
\ee
where the equality
$
V_rr=V_{\eta}\eta_rr
$
has been used.

{

From (4.20) and (4.21)  the functions  $r=\breve r(x,\eta)$ and $V= \breve V(x,\eta)$ can be found upon integration:

\be
\ln r=
C_7
\int
\fr 1
{R_1\sinh\eta}
d\eta
\ee

and

\be
\ln V=
-
\int\fr{\fr{1}{H^2} \sinh^2\eta + \fr{1}{P}- C_7}
{R_1\sinh\eta}
d\eta.
\ee

\ff

\nin
{\bf Assertion 4.2} {\it The integrals} (4.22)  {\it and} (4.23)
{\it yield the general solutions for the dependence of the functions   $r=\breve r(x,\eta)$ and $V= \breve V(x,\eta)$
on $\eta$
compatible with the {\it Total Set} of equations.
}

\s \s


\setcounter{equation}{0}

\s \s

\nin
{\bf \large 5  $\eta$-regular solution }

\s

\s

\nin
In the preceding Sect. 4 the equations have been derived systematically to represent the pseudo-Finsleroid metric
functions of the angle-separated type. In the present section, we find an interesting and significant  particular solution to the equations.
To specify the solution,  we introduce the following two definitions.

\s

\s

\nin
{\bf Definition 5.1} The  pseudo-Finslerian metric function  $F$ of the angle-separated type
is called {\it $\eta$-regular}, if
the characteristic functions
$V=\bvV(x,\eta)$ and $ r=\bvr(x,\eta)$
which enter the function
$F$
in accordance with the Separation Scheme (3.2)
are smooth of class $C^{\infty}$ with respect to their angle argument
$\eta$ over the total definition range $0 <\eta<\infty$.
The metric functions possessing such a property
will be denoted by $\FR$.

\s

 \s

\nin
{\bf Definition 5.2} If a metric function   $\FR$
 possesses the  property that the indicatrix is a space of constant negative curvature, then
the function is called the {\it $\eta$-regular pseudo-Finsleroid metric function},
to be denoted by  $\FRR$.

\s

\s

\nin
To make search for a  function $\FRR$,
let us choose  the {\it Class} II  among the three classes formulated below Formula (4.5) in Sect. 4.
Accordingly, we assume that $P>1$. Let us also denote
\be
\hC=\fr1{C},
\ee
where $C$ is the constant of
  Separation of Variables  that was arisen in the representation
$
\phi_t \phi_t=  C  Z_{tt}/Z
$
(see (4.9) in Sect. 4),
and introduce the constant $T$ by the help of the equality
\be
P=\fr1{T\hC}.
\ee

The appropriate idea is to subject the introduced constants to the inequalities
\be
T>1, ~ ~ 0<\hC<1, ~ ~ T\hC<1.
\ee
The inequality
$
H\ge 1
$
is always assumed.

{

We start with the representation
\be
F=b\bvV, \quad
\bvV=
C_1
\fr1
{R_1}
J,
\ee
and make  the choice
\be
C_7=\fr{1}{P},
\ee
where $R_1$ is the known function that was written out in (4.18) and $J=J(x,\eta)$ is
the input function to be determined.

\s

With the choice $C_7=1/P$, the equation  (4.21) reduces to merely
\be
\fr1{\bvV}
\bvV_{\eta}
=
-\fr{1}{H^2}
\fr1{R_1}
\sinh\eta
\ee
and
the equation  (4.20) takes on the  form
\be
\eta_rr=PR_1\sinh\eta.
\ee
The simple and useful equality
\be
\fr1{H^2}\lf({\eta}_r\rg)^2=
-\fr1VV_{rr}
\ee
can be derived.

{

From  the  equation (5.6), the dependence of the function $J=J(x,\eta)$ on $\eta$  can explicitly be found.
To this end we differentiate the equality  $\ln  \bvV=\ln(C_1)-\ln(R_1)+\ln J$ with respect to $\eta$ and compare the obtained result with (5.6). This method leads to the equation
\be
  (\ln J)_{\eta} =
  \lf(1-\fr1{H^2}\rg)
  \fr1{\sqrt{\hL}}   \sinh\eta,
\ee
where $\hL$ is the function given by (4.5).
 Integrating yields
\be
J=
 \exp\lf(
{H_1}\arcsinh \lf( \hL_1\cosh\eta\rg)
\rg),
\ee
where

\s

\be
{H_1}=\sqrt{1-\fr1{H^2}},  ~~  \hL_1=\fr {H_1}{S_1}, ~~   S_1=\sqrt{1-\fr1{P}}.
\ee
The representation (5.10) can alternatively  be written  in the form

\s

\be
J=
\lf(  \hL_1\cosh\eta+\sqrt{ (\hL_1)^2\cosh^2\eta+1}\rg)^{H_1}.
\ee
The right-hand part in (5.12) is obviously smooth of class
$C^{\infty}$ on the total range $0<\eta<\infty$.

\s

{

Next,  the function $r=\bvr(x,\eta)$ can be determined from (5.7). Indeed,
it is convenient to use the representation
\be
r=
\breve C_2
\fr {\sinh\eta}{R_1}
\,
Y_1, ~~  \breve C_2=\breve C_2(x).
\ee
It can readily be verified that  the introduced function $Y_1$ should be subjected to the equation

\s

\be
(\ln  Y_1)_{\eta}=
\lf( \fr1{P}-1\rg)
\fr1{\sinh\eta\sqrt{\hL}}.
\ee

\s

The attentive integration yields

\s

\be
Y_1=\lf(
\fr{\lf( 1-Nx+\sqrt{N+1}\sqrt{1+Nx^2}    \rg)(  x-1 )}
{\lf( 1+Nx+\sqrt{N+1}\sqrt{1+Nx^2}    \rg)(  x+1 )}
\rg)^{-\fr12\sqrt{1-\fr1{P}}    },
\ee
where
$ x=\cosh\eta $
and
\be
N=\fr{1-\fr1{H^2}}{\fr1{H^2}-\fr1{P}}.
\ee

{

From   (4.11) and (4.14) it follows that
\be
\bar\theta_f  f = \fr1{\hC}  R_2\sin\theta, ~~
\fr1U U_ff=   T\sin^2\theta,
\ee
 which in turn entails
\be
\fr 1{\bvU} \bvU_{\theta}=
\fr{1}{P}
\fr{   \sin\theta}
{
R_2}.
\ee

{

\s  \s

\nin
{\bf Definition 5.3} If a metric function   $\FRR$
 possesses the  property that the indicatrices in horizontal sections are of constant positive curvature, then
the function is called the {\it $\eta$-regular  Finsleroid-in-pseudo-Finsleroid metric function},
to be denoted by  $F^{\cal REG-FRD-IN-PSEUDO-FRD}$.

\s \s

\nin
In close resemblance with the function $R_1$ defined by (4.18),
we introduce the function
$$
R_2=
  \cos\theta
+
\sqrt{\hC}
\sqrt {  (T-1)\sin^2\theta+\lf(\fr1{\hC}-1\rg) \cos^2\theta },
$$
which can also be written  conveniently as
\be
R_2=
  \cos\theta
+
\sqrt { L_9 }, ~~ ~~  L_9= T\hC-\hC+(1-T\hC) \cos^2\theta .
\ee
Since  $ T>1$ and $T\hC<1$ (according to (5.3)),  the function $R_2=R_2(x,\theta)$ is totally regular with respect to the angle argument $\theta$,
namely is of  the class $C^{\infty}$.

From (5.17) the following simple formula can be obtained:
\be
U_{ff}=
(\theta_f)^2
T\hC U.
\ee

{

With this function $R_2$,   the dependence of the function $f=\bvf(x,\theta)$ on $\theta$
can be found from (5.17) with the help of the substitution
\be
f=
C_{17}
 \fr{\sin\theta}{R_2}
Y_2, ~~ C_{17}=C_{17}(x).
\ee
We obtain
\be
(\ln Y_2)_{\theta}=
-
(1-\hC)
\fr1{\sqrt{L_9}\sin\theta}.
\ee
The integration leads to the following result:

\s

\be
 Y_2= \lf(\fr{ \sqrt{\tan^2\theta+A}+\sqrt A}{ \sqrt{\tan^2\theta+A}-\sqrt A}\rg)^{\fr12\sqrt {1-\hC}},
 ~~~   A=\fr{1-\hC  }  {\hC T-\hC  }.
\ee

{

It is convenient to apply such a method  also to the function $U=\bvU(x,\theta)$.
Namely, postulating the representation
\be
U=C_{39}
\fr1{R_2}
I,  ~~  C_{39}=C_{39}(x),
\ee
from (5.18) we get
\be
(\ln I)_{\theta}=
 -(1-T\hC)
\fr{ \sin\theta}{\sqrt{L_9}}.
\ee
The integration shows that

\s

\be
I=
\lf(  \sqrt{1-T\hC}   \cos\theta
+
\sqrt {L_9}
 \rg)
^{ \sqrt{1-T\hC}}.
\ee

\ff

The interesting simple equality
\be
\lf(\ln \fr I{Y_2}\rg)_{\theta}=
\fr{\sqrt{L_9}}{\sin\theta}
\ee
can also be derived.

{

\s

It remains to consider the dependence of the function $\phi=\bar\phi(x,t)$ on the argument $t$,
where $t=w_1/w_2$ (see the definition of the variable $t$ in (3.3)).
According to (4.9), the function $\bar\phi$ must fulfill the equation
$$
\bar\phi_t \bar\phi_t=  \fr1{\hC}  \fr {Z_{tt}}Z,
$$
Let us solve the equation with the help of the function
\be
\phi=\fr1{\sqrt{\hC} }\arctan t +C^*, ~~ ~~   C^*=C^*(x),
\ee
which inverse is $ t=\tan(\sqrt{\hC}\phi-C^*) $,
so that
$
\bar\phi_t=  \lf(1/\sqrt{\hC} \rg) ( 1/(1+t^2)).
$
For the function $Z=Z(x,t)$ we obtain the equations
$$
 \fr1ZZ_{tt}=\fr 1{(1+t^2)^2}, ~~ \fr1ZZ_t=\fr t{1+t^2}\equiv \fr {w_1w_2}{w_1w_1+w_2w_2},
$$
which can be solved as follows:
\s

\be
Z=\sqrt{1+t^2}\,C_{11}\equiv \fr1{\cos(\sqrt{\hC}\phi-C^*) }\,C_{11}, ~~ ~~  C_{11}=C_{11}(x).
\ee
Recollecting the definition
$
f=w_2Z/w_3$ (see (3.2)),
we arrive at the representation
\be
 f=\fr1{w_3}w_{\perp} \,C_{11}(x),
\ee
where
\be
w_{\perp} =\sqrt{w_1w_1+w_2w_2}.
\ee
This  observation motivates introducing the following definition.
\s  \s

\nin
{\bf Definition 5.4} The Finslerian metric function $F$ is called $i_{\{3\}}$-axial  if
the the 1-forms $i=i_iy^i$ and $j=j_iy^i$ introduced in Sect. 1 enter the function  in the
 sum of squares $(i)^2+(j)^2$.

 \s \s

 \nin
   Thus, we are entitled to formulate  the following implication.

 \s \s

\nin
{\bf Assertion 5.1}
{\it
When the function $\phi=\phi(x,t)$ is given by Formula} (5.28),
{\it the function $f$ is $i_{\{3\}}$-axial, for  the variables $w_1,w_2$ enter the function in the
 sum of squares $(w_1)^2+(w_2)^2$.
}

\s \s

\nin
{\it Remark 5.1}
Formula (5.28) extends its pseudo-Riemannian precursor by presence  of the factor
$1/\sqrt{\hC}.$
The axial nature of the function $f$ is retained valid under the extension from the pseudo-Riemannian
geometry to the pseudo-Finsleroid theory under development in the present paper.
According to Formula (5.29), the form of dependence of the function $Z$ on the argument $t$ is exactly the same as in  the pseudo-Riemannian geometry.

\s  \s

{

\nin
If we now remind the definition $r=w_3U(x,f)$ (see (3.2)), we conclude that the function $r$ also is of the $i_{\{3\}}$-axial structure:
\be
r=r^*(x,w_3,w_{\perp}).
\ee
Because of the representation chain
$\{F=bV, V=\breve V(x,\eta), \eta=\eta(x,r), \theta=\bar\theta(x,f)\}$,
all the functions $F$, $V$,  $\eta$,  and $\theta$  are    $i_{\{3\}}$-axial.

{

\s

Differentiating the equality $ f=(1/w_3)w_{\perp} \,C_{11}$ with respect to $w_1, w_2$, and $w_3$
yields
\be
ff_1=\fr1{w_3w_3}w_1(C_{11})^2, ~~ ff_2=\fr1{w_3w_3}w_2(C_{11})^2,
~~ ff_3=-((w_1)^2+(w_2)^2)
\fr1{w_3w_3w_3}(C_{11})^2,
\ee
which entails
$
 f_1w_1+f_2w_2+f_3w_3=0.
$

\s


We can take $U_f$ from (5.17):
$
(1/U) U_ff=   T\sin^2\theta.
$
We get

$$
\D U{w_3}=U_ff_3
=-UT
\fr1f
\sin^2\theta f\fr1{w_3}
=-UT
\sin^2\theta \fr1{w_3},
$$
or
\be
\D U{w_3}=
-UT
\sin^2\theta \fr1{w_3}.
\ee
Here, the right-hand part does not involve neither  $R_2$  nor $I$.


Applying this Formula  to $r=w_3U$ yields

$$
r^*_{w_3}=U-
UT
\sin^2\theta , ~~    f=
C_{17}
 \fr{\sin\theta}{R_2}
Y_2,  ~~
  U=C_{39} \fr1  {  R_2}  I,
$$
from which it follows that

$$
r^*_{w_3}=
U- UT R_2R_2 f^2 \fr1{(C_{17})^2(Y_2)^2}
=
U- \fr1 U
\fr  {(w_{\perp})^2}{({w_3})^2}
\fr1{(Y_2)^2}
I^2 \,(C_{11})^2T
\fr{(C_{39})^2}{(C_{17})^2}.
$$
We obtain the representation

\be
r^*_{w_3}=
U- \fr1 U
\fr  {(w_{\perp})^2}{({w_3})^2}
\fr1{(Y_2)^2}
I^2 \,(C_{11})^2T
\fr{(C_{39})^2}{(C_{17})^2}.
\ee


Similarly, the derivative
$
\partial U/\partial {w_{\perp}}=C_{11} U_f/{w_3}
$
can readily transformed to the representation

\be
\D U{w_{\perp}}=UT
\sin^2\theta
\fr1{w_{\perp}}
\ee
which, when applied to
 $r=w_3U$, yields the result
$
r^*_{w_{\perp}}=
UT
(\sin^2\theta)
{w_3}/{w_{\perp}},
$
or
\be
r^*_{w_{\perp}}=
 \fr1 U
\fr  {w_{\perp}}{{w_3}}
\fr1{(Y_2)^2}
I^2 \,(C_{11})^2T
\fr{(C_{39})^2}{(C_{17})^2}.
\ee
The identity
$
w_3r^*_{w_3}+r^*_{w_{\perp}}w_{\perp}=r\equiv w_3U
$
is valid.

{

For tensorial evaluations, it is convenient to make the redefinition
$$
v^a=w_a
$$
(a=1,2,3) for  variables
and apply the representations

\be
r=\hat r(x,v),  ~~ f=\hat f(x,v),  ~~  U=\hat U(x,v),  ~~  \theta=\hat \theta(x,v),  ~~   \phi=\hat \phi(x,v), ~  ~~  v=\{v^a\}
\ee
to use  the derivatives
 \s

 $$
r_a=\D {\hat r}{v^a},   ~~ r_{ab}=\D {\hat r_a}{v^b}, ~~  \theta_a=\D {\hat \theta}{v^a}, ~~   \phi_a=\D {\hat \phi}{v^a},
~~
f_a=\D {\hat f}{v^a}, ~~ f_{ab}=\D {\hat f_a}{v^b}.
$$

Since  $r=v^3\hat U$, the function $\hat r$
is positively homogeneous of the first degree with respect to the variables $v^a$,
 and has the geometrical meaning of the {\it Finsler metric function in the horizontal sections}.
 The associated Finsler metric tensor is
\be
R_{ab}=rr_{ab}+r_ar_b
\ee
and $\{r_a\}=\{r_1,r_2,r_3\}$ plays the role of the respective covariant unit vector; $r_av^a=r$.

\s

Because of the axial structure (5.32) of $r$,  the function $\bar r(x,v^3,v_{\perp})$ exists
such that $r=\bar r$ and
$$
\hat r=\bar r(x,v^3,v_{\perp}),
$$
where
$$
v_{\perp}=\sqrt{  (v^1)^2+(v^2)^2 }
$$
(which is equal to $w_{\perp}$).
We derivatives
$$
r_1=\bar r_{ v_{\perp} }\fr{v^1}{ v_{\perp} }, ~~  r_2=\bar r_{ v_{\perp} }\fr{v^2}{ v_{\perp} }, ~~
r_3=\bar r_{ v^3}
$$
are convenient to use in calculations.

 {

We can readily obtain
$$
f_1f_1+ff_{11}=\fr1{v^3v^3}(C_{11})^2, ~~  f_1f_2+ff_{12}=0, ~~
 f_1f_3+ff_{13}=
-2\fr1{v^3v^3v^3}v^1(C_{11})^2,
$$

\s

$$
f_2f_1+ff_{21}=0, ~~  f_2f_2+ff_{22}=\fr1{v^3v^3}(C_{11})^2, ~~
 f_2f_3+ff_{23}=
-2\fr1{v^3v^3v^3}v^2(C_{11})^2,
$$
and

\s

$$
f_1f_3+ff_{13}=-2v^1 \fr1{v^3v^3v^3}(C_{11})^2, ~~  f_2f_3+ff_{23}=-2v^2 \fr1{v^3v^3v^3}(C_{11})^2,
$$
together with

\s

\be
 f_3f_3+ff_{33}=3\fr{((v^1)^2+(v^2)^2)}{v^3v^3v^3v^3}(C_{11})^2.
\ee

{

The following assertion is valid.

\s

\nin
{\bf Assertion 5.2} {\it Given  the function $R_2$ by the help of   Formula} (5.19).
{\it If the   equations
\be
\fr 1U U_{\theta}=
\fr{1}{P}
\fr{   \sin\theta}
{
R_2},
 ~ ~
f_{\theta}
=
\hC
f
\fr1{\sin\theta R_2}
\ee
are fulfilled, then the indicatrices of the sectional horizontal spaces
possess the property of positive curvature, and {\it vice versa}.
The respective curvature value $\cC^{\{3\}}(x)$ will be indicated below in Formula }(5.56).

\s  \s


{

To verify  the assertion we can evaluate the tensor
$
C_{abc}=(1/2)\partial{R_{ab}}/\partial{v^c}.
$
The insertion of the  metric tensor
$
R_{ab}=r_ar_b+rr_{ab}
$
(see (5.39))
yields the symmetric   representation
\be
2C_{abc}=r_ar_{bc}+r_br_{ac}+r_cr_{ab}  +   rr_{abc}.
\ee
Raising the indices is naturally performed by the help of
the tensor $R^{fd}$  reciprocal to $R_{fd}$, as exemplified by
$C^f_{be}=R^{fd}C_{dbe}$.
In the process of evaluating  the curvature tensor
$$
  R^*_{bace} = C_{afc}C^f_{be}-C_{afe}C^f_{bc},
$$
it is appropriate to  apply   Formulas (5.34)-(5.40) and their concomitants,
and  in the last steps   use the equalities
$
U_{ff}= (\theta_f)^2T\hC U
$
(see (5.20)),
$
(1/H^2)\lf({\eta}_r\rg)^2=  -(1/V)V_{rr}
$
(see (5.8)), and
$$
\sin\theta
\fr1{(\theta_f)^2}\lf(\theta_{ff}+2T\sin^2\theta \fr1{f}\theta_f\rg)=
\sqrt{L_9}
-\fr{(1-T\hC)  \sin^2\theta}{\sqrt{L_9}}
$$
(this equality can be verified with the help of (5.19)).
This process (which is not short) leads
to the remarkable representation
\be
r^2   R^*_{bace}=(P-1)(h_{bc} h_{ae} - h_{be} h_{ac}),
\ee
where $h_{bc}=R_{bc}-r_br_c$.

{

Also, we can obtain
the angler representation
\be
h_{ab}=
\fr {1}P
(\theta_a\theta_b+\sin^2\theta\phi_a\phi_b)
r^2
\ee
by following closely the  calculation method  disclosed and applied in the previous paper [1,2].

\ff

{

Let us apply these observations to elucidating the curvature properties of the  horizontal sections.

Like to Sect. 3 in [1,2], we denote by
\be
{\cal T}_{x;\la(x)}\in T_xM
\ee
 the space orthogonal to the vertical $b$-axis
  at a fixed value $\la_x=b_x.$
This 3-dimensional space ${\cal T}_{x;\la(x)}$
is the  horizontal section of the tangent space $T_xM $ through  the axis point assigned by
a value $\la(x)=\la_x$.

\s

The sectional horizontal space ${\cal T}_{x;\la(x)}$ consists of the end points of the tangent vectors
$y_x\in T_x$ representable by the component expansion
\be
y_x^i=\la(x) b^i(x)+ii^i(x)+jj^i(x)+ i_{\{3\}}  i_{\{3\}}^i(x),
\ee
where
$\la(x), b^i(x), i^i(x), j^i(x),   i_{\{3\}}^i(x)$
are  fixed and  $i,j,i_{\{3\}}$ are arbitrary.
We can identify
$b_x=b_i(x)y_x^i.$

Following the geometrical imagination adopted in Sect. 3 of [1,2],
we introduce the notion of
the 3-dimensional { \it centered vector  space}
\be
{\cal  V}^{\{3\}}_{x;\la(x)}=\{O_{x;\la(x)}, {\cal T}_{x;\la(x)}, ~~ u\in {{\cal T}_{x;\la(x)}},  ~   u=\{u^a\}, ~ a=1,2,3\},
\ee
which is of the  sectional horizontal meaning, where
$
u^1=i, ~  u^2=j, ~  u^3=i_{\{3\}},
$
and    $\{u^a\}$  are geometrically interpreted as the components of the  vectors  $u\in {{\cal T}_{x;\la(x)}}$
that are supported by  the center point $O_{x;\la(x)}$ of the space  ${\cal  V}^{\{3\}}_{x;\la(x)}$. The entered $O_{x;\la(x)}$ belongs to the vertical $b$-axis,
namely  $O_{x;\la(x)}$ is the point at which ${\cal T}_{x;\la(x)}$ intersects this axis.

{

The function
$
r=\hat r(x,v)
$
introduced in (5.38) determines the function

\be
F_{x;\la(x)}^{\{3\}}=F_{\la(x)}^{\{3\}}(x,u) ~~ \text {with} ~~ F_{\la(x)}^{\{3\}}(x,u)=\hat r_{\la(x)}(x,u)
\ee
which we shall naturally interpret as  the {\it Finsler metric function in the sectional horizontal space}
${\cal T}_{x;\la(x)}.$

If we take into account the first-order homogeneity of the function $\hat r(x,v)$ with respect to $v=\{v^a\}$
(which is obvious from the input definition  $r=v^3U$ (see (3.2)) for the function $r$), together with the relation $v^a=u^a/\la$,
we establish the important equality
\be
\la(x)\hat r_{\la(x)}(x,v)=\hat r_{\la(x)}(x,u).
\ee

{

In the tangent space $T_x$   we confine  our analysis to
region
bounded by the indicatrix
$\cI_x$  (introduced by Definition 1.1 in  Sect. 1).
Accordingly,
  we should restrict  the consideration
 by the $b$-like parts
\be
{\cal PT}_{x;\la(x)} = \{y\in {\cal PT}_{x;\la(x)}: y\in {\cal T}_{x;\la(x)}, ~ F_{x;\la(x)}^{\{3\}}\le {\cal R}_{x;\la(x)}\}
\ee
of the horizontal sections ${\cal T}_{x;\la(x)}$, where
\be
{\cal R}_{x;\la(x)}=\lf(F_{x;\la(x)}^{\{3\}}(x,u)\rg)_{x;\la(x);max}
\ee
is the {\it radius of the horizontal section} ${\cal PT}_{x;\la(x)}$.

The maximal value
$
\lf(\hat r(x,v)\rg)_{\la(x);max}
$
of the function $\hat r(x,v)$ in ${\cal PT}_{x;\la(x)}$ obeys obviously  the equation
$$
\la(x)V\lf(x,\lf(\hat r(x,v)\rg)_{\la(x);max}\rg)=1.
$$
With this value, the radius can be  given by the equality
\be
{\cal R}_{\la(x)}=\la(x)\lf(\hat r(x,v)\rg)_{\la(x);max},
\ee
because
$$
\lf( F_{x;\la(x)}^{\{3\}}(x,u)\rg)_{max}=\lf(\hat r(x,u)\rg)_{\la(x);max}
$$
and the equality (5.49) holds identically.
\s

Accordingly, we  restrict the centered vector space  ${\cal V}^{\{3\}}_{x;\la(x)}$ by the part
${\cal P V}^{\{3\}}_{x;\la(x)}$ which is bounded by
 the indicatrix
$\cI_x$.
The space ${\cal P V}^{\{3\}}_{\la(x)}$ is obtained by substituting
 ${\cal PT}_{x;\la(x)}$
 with   ${\cal T}_{x;\la(x)} $ in the definition (5.47)
 of ${\cal V}^{\{3\}}_{x;\la(x)}$.

{

From  the function $F_{x;\la(x)}^{\{3\}}=F_{x;\la(x)}^{\{3\}}(x,u)$ we construct the Finsler metric tensor
$G_{x;\la(x);ab}=G_{x;\la(x);ab}(x,u)$ according to the ordinary Finslerian rule
\be
G_{x;\la(x);ab}=\fr12 \Dd {\lf(F_{x;\la(x)}^{\{3\}}\rg)^2}{u^a}{u^b}.
\ee
The    sectional horizontal space ${\cal P V}^{\{3\}}_{x;\la(x)}$ becomes
the three-dimensional   {\it sectional   Finsler space}
 \be
{\cal  F}^{\{3\}}_{x;\la(x)}=\{{\cal P V}^{\{3\}}_{x;\la(x)},~ G_{x;\la(x);ab},  ~  F_{x;\la(x)}^{\{3\}}  \le {\cal R}_{x;\la(x)}\}.
\ee
The  axis of this space is assigned by
the direction of the   vector  $i^i_{\{3\}}(x)$.

  \s  \s

\nin
 {\bf  Definition 5.5}
The surface $   \cI_{x;\la(x)}^{\{2\}}   \subset{\cal P V}^{\{3\}}_{x;\la(x)}$ introduced by
 \be
\cI_{x;\la(x)}^{\{2\}} = \{  u \in  {\cal P V}^{\{3\}}_{x;\la(x)}, ~  F_{\la(x)}^{\{3\}} (x,u) = {\cal R}_{x;\la(x)}   \}
\ee
is called the {\it indicatrix of the sectional Finsler space}.

\s  \s

\nin
${\cal R}_{x;\la(x)}   $ is the radius of this indicatrix,
which  is the boundary  of the centered space ${\cal P V}^{\{3\}}_{x;\la(x)}$.

{

Considering the equalities
$$
G_{x;\la(x);ab}=\fr12 \Dd {\lf(F_{x;\la(x)}^{\{3\}}\rg)^2}{u^a}{u^b}
=\fr12 \Dd {\lf( \hat r_{\la(x)}(x,u)\rg)^2}{u^a}{u^b}
=\fr1{2} \Dd {\lf( \hat r_{\la(x)}(x,v)\rg)^2}{v^a}{v^b},
$$
where (5.49) has been used,
we arrive at  the identification
$$
G_{x;\la(x);ab}=
 R_{ab},
$$
in which the last tensor $R_{ab}$ is known from (5.39). With the knowledge of the tensor,
and taking into account Formulas (5.43) and (5.44),
the subsequent attentive evaluation leads to

\be
\cC_{x;\la(x)}^{\{3\}}=\fr{P(x)}{\lf({\cal R}_{x;\la(x)}\rg)^2}.
\ee
Assertion 5.2 is valid.
The obtained  formula (5.56) complies with the well-known geometrical property that the curvature of the Euclidean sphere of radius $R$ equals $1/R^2$.

{

The evaluation of the determinant of the tensor $R_{ab}$ yields the result

\be
\det(R_{ab})=
\fr{1}{P^2\hC^2}
I^6
\fr1{(Y_2)^4}
\lf(\fr{C_{11}}{C_{17} }\rg)^4
(C_{39})^6>0
\ee
which is positive and of the angle-regular  class $C^{\infty}$ with respect to the involved angles
$\eta$ and $\theta$, and is independent of the  angle $\phi$ because the space is  $i_{\{3\}}$-axial.

{

Thus we are justified to formulate the following statements.

\s  \s

\nin
{\bf Assertion 5.3}
{\it If the characteristic functions $\{V=\breve V(x,\eta),r=\breve r(x,\eta),U=\breve U(x,\theta),f=\breve f(x,\theta)\}$ in
a   pseudo-Finslerian metric function $F(x,y$) of the angle-separated  $i_{\{3\}}$-axial type
fulfill the ordinary differential equations
\be
\fr1V
V_{\eta}=
-
\fr1{H^2}
\fr1{ R_1}
\sinh\eta,
 ~ ~
r_ {\eta}
=\fr rP
 \fr 1{\sinh\eta R_1},
 ~ ~
\fr 1U U_{\theta}=
\fr{1}{P}
\fr{   \sin\theta}
{
R_2},
 ~ ~
f_{\theta}
=
\hC
f
\fr1{\sin\theta R_2},
\ee
where the subscripts  $\eta$ and $\theta$ mean the differentiations,
then:

\s

\nin
the function $F$ is  of the  $ \FRRIN$-type;

\s

\nin
the indicatrix  $\cI_x$ is of the constant negative curvature which value is $-H^2$;

\s

\nin
in the horizontal sections,
the geometry is assigned by  the Finslerian metric function ${\cal  F}^{\{3\}}_{x;\la(x)}$,
  the  associated indicatrix  $\cI_{x;\la(x)}^{\{3\}}$  is  of the constant positive  curvature
 $\cC_{x;\la(x)}^{\{3\}}$ given by} (5.56),
{\it
and the metric tensor is angle-regular and positive-definite. }

\s  \s

{

\s  \s

\nin
{\it Remark 5.2} The class $C^{\infty}$ regularity of the dependence of the function $\bvV$ on $\eta$
over all the definition range $\eta\in (0,\infty)$ is evidenced from the right-hand part
of the formula
$(1/{\bvV})\bvV_{\eta}=-1/(H^2R_1) \sinh\eta $
 displayed in (5.6), because the  function $R_1$  obviously possesses such a regularity property
(see the representation for $R_1$ in (4.18)). The representation (5.10) of the function $J$ also manifests
the property. The regularity formulated comes actually from the perfect regularity of the right-hand part
in the
input representation (4.5) of the function $\hL>0$, because we are keeping the inequalities $P>1$ and $H>1$ well.
The similar motivation can be addressed to the function
$\bvU=\bvU(x,\theta)$,  namely we can refer to the formula
$(1/{\bvU})\bvU_{\theta}=1/(PR_2) \sin\theta $ shown in (5.18)
and take into account the representation (5.19) which introduces the function $R_2$ with the help of the function
$L_{9}$. Upon the restrictions $ \{T>1, T\hC<1\}$ assumed the function
$L_{9}$  is positive and  also smooth of  class $C^{\infty}$ with respect to the angle argument $\theta$.

\s  \s


\nin
{\bf\large Conclusions: comparisons and properties}

\s  \s

\nin
Thus, we have at our disposal two metric functions of the Finsleroid-in-pseudo-Finsleroid-type,
namely $F^{\{P<1\}}$ (Class I) and $F^{\{P>1\}}$ (Class II).
By $F^{\{P<1\}}$ we denote the function proposed and described in the previous paper [1] (which was preceded
by the publication [2]).
$F^{\{P>1\}}$ is the function that was found and investigated in Sec. 5 of the present paper.

The angle separation  (3.1)-(3.5) is equally meaningful  for both the metric functions.
The expansion
$$
l^i=\lf(b^i+w_1i^i+w_2j^i +w_3i_{\{3\}}^i\rg)\fr1V
$$
can universally be applied with
$$
w_3= r\fr1U, ~~ w_1=w_{\perp}\cos\phi,   ~~ w_2=w_{\perp}\sin\phi,
~~
 w_{\perp} =\fr1{C_{11}}  w_3f,
$$
where  (5.30) and (5.31) have been used.

Outwardly, the distinction of the  derivative  list
(5.58) of Assertion 5.3 from the case developed in the previous paper [1]
is only in the presence of the scalar $\hC=\hC(x)$ in the last member which represents
$f_{\theta}$ (see the list of Formulas (22)-(23) in [1]).
However, in each tangent space $T_x$ the constant  $\hC$
gives rise to essential changes in the characteristic functions, including
 the appearance of new function   $Y_2$ which enters Formula (5.21) for $f=f(x,\theta)$.
The dependence of $Y_2$ on $\theta$ is rather complicated (see (5.23)), which does not permit to
obtain the inverse function $\theta=\theta(x,f)$ in an explicit algebraic form
(which was possible in the paper [1],  because of the presumption  $\hC=1$ made in [1]).


The significant distinction between the functions $F^{\{P>1\}}$ and $F^{\{P<1\}}$ is rooted in the structure of the functions
$\hL$ and $L$ which enter the key function $R_1$, namely
in the present paper we used
$$
 R_1= \cosh\eta  +\sqrt{\hL},  ~~
\hL=1-\fr1{P}+\lf(1-\fr1{H^2}\rg)  \sinh^2\eta , ~~ P>1
$$
(see Formulas (4.5) and (4.18) in Sect. 4), while
in  [1]  the function $R_1$ was taken  according to
$$
 R_1= \cosh\eta  +\sqrt{L}, ~~
 L=1-\fr1{p^2}+\lf(1-\fr1{H^2}\rg)  \sinh^2\eta , ~~ p^2<1
$$
(Formula (12) in [1]).
Indeed, the involved function $\hL$, and whence the function $R_1$, is totally positive and regular over
$\eta\in (0,\infty)$. This observation explains why the pseudo-Finsleroid metric function proposed in Sect. 5 is
$\eta$-regular.

Let us substitute  the notation $P$ with $p^2$ in the second version of $R_1$.
Since in this case $P<1$, the function $L$ vanishes at the value $\eta_0>0$ obtainable from the algebraic equation
$$
\fr1P-1=\lf(1-\fr1{H^2}\rg)  \sinh^2\eta_0,
$$
so that $L<0$ if $\eta<\eta_0$, and $L>0$ if $\eta>\eta_0$. Negative $L$ are
nonadmitted  under the square root in the right-hand part of $R_1$, whence we encounter with
the
{\it deficit of angle} $\eta$:  the region $\eta<\eta_0$
is inaccessible.
The internal cone appears, because the function $F^{\{P<1\}}$  vanishes at $\eta=\eta_0$,
so that $F^{\{P<1\}}$ is the  two-cone pseudo-Finsleroid metric function (see [1]).
In case of $F^{\{P>1\}}$, only one cone is evidenced.


In the process of the  systematic evaluation performed in Sect. 5 to obtain the metric function
 $F^{\{P>1\}}$,
 we have found explicitly the dependence of all the involved functions $V,r,U,f,Z$ on the angle triple $\eta,\theta,\phi$. The knowledge of this dependence, when taken in conjunction with the representation
 of the unit vector $l^i$ which has been written out in the beginning of the present Section,
yields the explicit dependence $
y^i=\lf(b^i+w_1i^i+w_2j^i +w_3i_{\{3\}}^i\rg)b
$
 of tangent vectors $\{y^i\}$ on the angle triple.
At the same time,
  the inverse explicit dependence
  can be found  for the function $\phi=\phi(x,y)$ (see (5.28)) and cannot be obtained for the functions
  $\eta=\eta(x,y)$ and $\theta=\theta(x,y)$.
 The generating  function $V$ which enters the product  $F=bV$  is specified by the representation
(5.4) which is of the type $V=V(x,\eta)$. Therefore, the function $V$ borrows from $\eta$ the implicit character of dependence on tangent vectors $y$.
Nevertheless, the derivatives
$\partial V/\partial \eta$, $\partial r/\partial \eta$, $\partial U/\partial \theta$, $\partial f/\partial \theta$,
as well as
 $\partial V/\partial y^i$, $\partial\eta/\partial y^i$,  $\partial\theta/\partial y^i$, and $\partial U/\partial y^i$   admit simple explicit algebraic representations involving the angle triple  $\eta, \theta,\phi$
  (see the list (5.58) in Assertion 5.3).
  The latter property opens up a direct way to evaluate  the components of the  covariant unit vector
  $l_i=\partial F\partial y^i$
  and metric tensor $\{g_{ij}\}$ in  concise  forms. With these key objects at hands, we are able to clarify the structure of  the indicatrix metric tensor and many other objects.

The (Class II)-metric function $F^{\{P>1\}}$
shears with  (Class I)-metric function $F^{\{P<1\}}$ the
smoothness  of class $C^2$
on all the subspace
  ${\mathbf{\cal IF}}_{ \{x\}} \subset T_xM$
 bounded by the pseudo-Finsleroid.
Violations of differentiability  meet only on
the  two-axes section   ${\cal S}_{b,i_{\{3\}}{ \{x\}}  }\subset  T_xM$.
The smoothness  of class   $C^{\infty}$ holds perfect on the space
${\mathbf{\cal IF}}_{ \{x\}} \setminus {\cal S}_{b,i_{\{3\}}{ \{x\}}  }$.
 If a  three-dimensional hyperplane in  the tangent space $ T_xM$
 includes the vertical $b$-axis then the hyperplane is called
the {\it vertical section} of $ T_xM$.
The  two-axes section ${\cal S}_{b,i_{\{3\}}{ \{x\}}  }\subset  T_xM$
 is the  vertical section of  $ T_xM$
if  the horizontal
 $ i_{\{3\}}$-axis also belongs
to the section (see more detail in [1,2]).

\s  \s


\setcounter{equation}{0}

\nin
{\bf \large A Appendix: Skew-Vanishing Lists}

\s \s

\nin
The  Skew-Vanishing Lists   $\cS_1, \cS_2,\cS_3$ displayed in Sect. 3
are implications from three integrability conditions entered (3.12).
Let us develop the first condition
\be
\D{\eta_{ij}}{y^n}-\D{\eta_{in}}{y^j}=0,
\ee
taking into account
the angle representation
$$
h_{ij}=
-
\fr{F^2}{H^2}
\Bigl(
\eta_i
\,\eta_j
+
\sinh^2\eta (
\theta_i\theta_j  +  \sin^2\theta  \phi_i\phi_j
)
\Bigr) \equiv F\D{l_i}{y^j}
$$
indicated in Formula (2.2) of Sect. 2.
Using the Particular Expansion (3.8) of $\eta_{ij}$, from (A.1)
we straightforwardly obtain

{

$$
0=\D{\eta_{ij}}{y^n}-\D{\eta_{in}}{y^j}=
\fr1{F^2}l_i(l_n\eta_j-l_j\eta_n)
-\fr1{F^2}(h_{in}\eta_j-h_{ij}\eta_n)
-\fr1F(l_j\eta_{in}-l_n\eta_{ij})
  $$

  \s \s

$$
+
2u_2\fr{1}{H^2}
\Bigl(
\sinh\eta \sin^2\theta\cosh\eta\eta_n +\sinh^2\eta \sin\theta\cos\theta\theta_n
\Bigr)
\phi_i\phi_j
$$

  \s \s

$$
-
2u_2\fr{1}{H^2}
\Bigl(
\sinh\eta \sin^2\theta\cosh\eta\eta_j +\sinh^2\eta \sin\theta\cos\theta\theta_j
\Bigr)
\phi_i\phi_n
$$

  \s \s

$$
+
u_2\fr{1}{H^2}
\sinh^2\eta \sin^2\theta
\phi_j
\Biggl[
-\fr1F(l_i\phi_n+l_n\phi_i)  +r_1\fr{1}{H^2}\sinh^2\eta\sin\theta(\phi_i\theta_n+\phi_n\theta_i)
 +
r_2\fr{1}{H^2}
\sinh^2\eta\sin^2\theta  \phi_i\phi_n
$$

\s \s

$$
+
r_5\fr{1}{H^2}
\sinh\eta\sin\theta (\phi_i\eta_n+\phi_n\eta_i)
 \Biggr]
$$

\s  \s

$$
-
u_2\fr{1}{H^2}
\sinh^2\eta \sin^2\theta
\phi_n
\Biggl[
-\fr1F(l_i\phi_j+l_j\phi_i)
+
r_1\fr{1}{H^2}\sinh^2\eta\sin\theta(\phi_i\theta_j+\phi_j\theta_i)
+
r_2\fr{1}{H^2}
\sinh^2\eta\sin^2\theta  \phi_i\phi_j
$$

\s \s

$$
+
r_5\fr{1}{H^2}
\sinh\eta\sin\theta (\phi_i\eta_j+\phi_j\eta_i)
 \Biggr]
$$

  \s  \s

{

$$
  +
  2 u_3\fr{1}{H^2}\sinh\eta\cosh\eta\theta_i
  (\theta_j\eta_n-\theta_n\eta_j)
$$

  \s \s

$$
  + u_3\fr{1}{H^2}\sinh^2\eta
\theta_j
\Biggl[
-\fr1F(l_i\theta_n+l_n\theta_i)
  +
 z_2\fr{1}{H^2}
\sinh^2\eta\sin^2\theta  \phi_i\phi_n
+
z_4\fr{1}{H^2}
\sinh\eta (\eta_i\theta_n+\eta_n\theta_i)
\Biggr]
$$

  \s  \s

$$
-
 u_3\fr{1}{H^2}\sin^2\eta
\theta_n
\Biggl[
-\fr1F(l_i\theta_j+l_j\theta_i)
  +
 z_2\fr{1}{H^2}
\sinh^2\eta\sin^2\theta  \phi_i\phi_j
+
z_4\fr{1}{H^2}
\sinh\eta (\eta_i\theta_j+\eta_j\theta_i)
\Biggr]
$$

  \s \s

{

$$
+
u_6\fr{1}{H^2}
\eta_j
\Biggl[
-
\fr1F(l_i\eta_n+l_n\eta_i)
+
u_2\fr{1}{H^2}
\sinh^2\eta \sin^2\theta \phi_i\phi_n  + u_3\fr{1}{H^2}\sinh^2\eta\theta_i\theta_n
+
u_6\fr{1}{H^2}\eta_i\eta_n
\Biggr]
$$

\s \s

$$
-
u_6\fr{1}{H^2}
\eta_n
\Biggl[
-
\fr1F(l_i\eta_j+l_j\eta_i)
+
u_2\fr{1}{H^2}
\sinh^2\eta \sin^2\theta \phi_i\phi_j
  + u_3\fr{1}{H^2}\sinh^2\eta\theta_i\theta_j
+
u_6\fr{1}{H^2}\eta_i\eta_j
\Biggr]
$$

\s  \s

$$
+
u_{2n}\fr{1}{H^2}
\sin^2\eta \sinh^2\theta \phi_i\phi_j
+
 u_{3n}\fr{1}{H^2}\sinh^2\eta\theta_i\theta_j
+
u_{6n}\fr{1}{H^2}\eta_i\eta_j
$$

  \s  \s

$$
-
u_{2j}\fr{1}{H^2}
\sinh^2\eta \sin^2\theta \phi_i\phi_j
-
 u_{3j}\fr{1}{H^2}\sinh^2\eta\theta_i\theta_n
-
u_{6j}\fr{1}{H^2}\eta_i\eta_n.
$$
After convenient rearrangements we get

{

\s  \s

$$
\fr1{F^2}l_i(l_n\eta_j-l_j\eta_n)
-\fr1{F^2}(h_{in}\eta_j-h_{ij}\eta_n)
-\fr1F(l_j\eta_{in}-l_n\eta_{ij})
  $$

  \s \s

$$
+
u_2
\lf(
-
u_6
\fr{1}{H^4}
\sinh^2\eta \sin^2\theta
+
2\fr{1}{H^2}
\sinh\eta \sin^2\theta\cosh\eta
+
r_5
\fr{1}{H^4}
\sinh^3\eta \sin^3\theta
\rg)
\phi_i(\eta_n \phi_j-\eta_j \phi_n)
$$

  \s  \s

$$
+
u_2
\lf(
2\fr{1}{H^2}
\sinh^2\eta \sin\theta\cos\theta
+
r_1\fr{1}{H^4}\sinh^4\eta\sin^3\theta
 \rg)
\phi_i(\theta_n\phi_j-\theta_j\phi_n)
$$

  \s \s

$$
+
u_2\fr{1}{H^2}
\sinh^2\eta \sin^2\theta
\fr1F\phi_i(\phi_nl_j-\phi_jl_n)
$$

  \s \s

$$
  +
 u_3
 \lf(
-
u_6\fr{1}{H^4}
\sinh^2\eta
 +
   2 \fr{1}{H^2}\sinh\eta\cosh\eta
+
z_4
\fr{1}{H^4}
\sinh^3\eta
\rg)
\theta_i
  (\theta_j\eta_n-\theta_n\eta_j)
$$

  \s  \s

$$
-
 u_3\fr{1}{H^4}
 z_2 \sinh^4\eta
\sin^2\theta
\phi_i( \theta_n  \phi_j-\theta_j  \phi_n)
+
 u_3\fr{1}{H^2}\sinh^2\eta
\fr1F
\theta_i(\theta_nl_j-\theta_jl_n)
$$

\s    \s

$$
+
u_6\fr{1}{H^2}
\fr1F
\eta_i(\eta_nl_j-\eta_jl_n)
$$

\s \s

$$
+
u_{2n}\fr{1}{H^2}
\sinh^2\eta \sin^2\theta \phi_i\phi_j
+
 u_{3n}\fr{1}{H^2}\sinh^2\eta\theta_i\theta_j
+
u_{6n}\fr{1}{H^2}\eta_i\eta_j
$$

  \s  \s

$$
-
u_{2j}\fr{1}{H^2}
\sinh^2\eta \sin^2\theta \phi_i\phi_j
-
 u_{3j}\fr{1}{H^2}\sinh^2\eta\theta_i\theta_n
-
u_{6j}\fr{1}{H^2}\eta_i\eta_n
=0.
$$

\s  \s

Here we apply the Symmetrizing Condition (3.34)
$z_2
\sinh^2\eta
\sin\theta
=
r_1
\sinh^2\eta \sin^2\theta
+
2H^2
 \cos\theta,
$
 such that we are left with

{

\s  \s

$$
u_2
\lf(
-
u_6
\fr{1}{H^4}
\sinh\eta
+
2\fr{1}{H^2}
\cosh\eta
+
r_5
\fr{1}{H^4}
\sinh^2\eta \sin\theta
\rg)
\sinh\eta \sin^2\theta
\phi_i(\eta_n \phi_j-\eta_j \phi_n)
$$

  \s  \s

$$
  +
 u_3
 \lf(
-
u_6\fr{1}{H^4}
\sinh\eta
 +
   2 \fr{1}{H^2}\cosh\eta
+
z_4
\fr{1}{H^4}
\sinh^2\eta
\rg)
\sinh\eta
\theta_i
  (\theta_j\eta_n-\theta_n\eta_j)
$$

  \s  \s

$$
+
u_{2n}\fr{1}{H^2}
\sinh^2\eta \sin^2\theta \phi_i\phi_j
+
 u_{3n}\fr{1}{H^2}\sinh^2\eta\theta_i\theta_j
+
u_{6n}\fr{1}{H^2}\eta_i\eta_j
$$

  \s  \s

$$
-
u_{2j}\fr{1}{H^2}
\sinh^2\eta \sin^2\theta \phi_i\phi_j
-
 u_{3j}\fr{1}{H^2}\sinh^2\eta\theta_i\theta_n
-
u_{6j}\fr{1}{H^2}\eta_i\eta_n
=
\fr1{F^2}(h_{in}\eta_j-h_{ij}\eta_n).
$$

{

Inserting

\s

$$
u_{2n}=u_{2\eta}\eta_n+u_{2\theta}\theta_n+u_{2\phi}\phi_n,  \s
u_{3n}=u_{3\eta}\eta_n+u_{3\theta}\theta_n+u_{3\phi}\phi_n,  \s
u_{6n}=u_{6\eta}\eta_n+u_{6\theta}\theta_n+u_{6\phi}\phi_n,
$$
reduces the equality under study to

  \s

$$
u_{2\theta}
=
u_{2\phi}
=
u_{6\theta}
=
u_{6\phi}
=0,
$$

\s  \s

$$
\lf(u_{2\eta}+\fr{1}{H^2}u_2(u_2-u_6)\rg)
\sinh^2\eta \sin^2\theta
\phi_i(\eta_n \phi_j-\eta_j \phi_n),
$$
and

\s

$$
  +
\lf(u_{2\eta}+\fr{1}{H^2} u_2(u_2-u_6)\rg)
\sinh^2\eta
\theta_i
  (\eta_n\theta_j-\eta_j\theta_n)
=
\fr{H^2}{F^2}(h_{in}\eta_j-h_{ij}\eta_n).
$$
The tensor $  h_{ij} $ can be expanded in accordance with
the angle representation indicated below (A.1).
In this way, the integrability condition (A.1) is reduced to the following set
\be
u_{2\theta}
=
u_{2\phi}
=
u_{6\theta}
=
u_{6\phi}
=0,  ~~
u_{2\eta}+\fr{1}{H^2}u_2(u_2-u_6)=1,
\ee
and
\be
z_{4\theta}
=
z_{4\phi}
=0, ~~
(r_5\sin\theta)_{\theta}
=
r_{5\phi}
=0.
\ee
Thus the List $\cS_1$  (see (3.13)-(3.15)) is valid.

The similar method of evaluation leads to verify the Lists $\cS_2$ and $\cS_3$.


\s \s

\nin
{\bf \large B Appendix: Evaluation of the  curvature tensor }

\setcounter{equation}{0}

\s  \s

\nin
Below we show how the equalities
\be
L_2=L_3, ~~ F^2(C^i{}_{pq}C_{ijn}-C^i{}_{pn}C_{ijq})
 =
T^*
(h_{pq}h_{jn}-h_{pn}h_{jq}),
\ee

\s
\be
T^*=-\fr1{H^2}
\Bigl( (L_2)^2
+
L_2
(u_6-2L_2)
\Bigr),
\ee
and
\be
T^*=1-H^2, ~~  L_2
(u_6-2L_2)
=
L(z_3\sinh\eta-L)
\ee
(see (3.39)-(3.44))
can be arrived at.

To this end we use  the expansion (3.36) of the tensor $C_{ijn}$ written  in  Sect. 3
in  terms of the unit vectors
and evaluate the contraction
$$
H^2F^2C^i{}_{pq}C_{ijn}=
\Biggl[
  u_6u^iu_pu_q
+
L_3
 (u^im_p m_q+u_qm_pm^i  +  u_pm^im_q)
$$

\s

\s

$$
+
L_2
(u^ip_p p_q  + u_qp^i p_p+u_pp^ip_q)
+
z_3
\sinh\eta
 m^im_pm_q
$$

\s

\s

$$
+
L
 (m_pp^ip_q+m_q p^ip_p+m^i p_pp_q)
+
 r_2
\sinh\eta
\sin\theta
p^ip_pp_q
\Biggr]
    \times
 $$

\s

\s

$$
\Biggl[
  u_6u_iu_ju_n
+
L_3
 (u_im_j m_n+u_nm_jm_i  +  u_jm_im_n)
$$

\s

\s

$$
+
L_2
(u_ip_j p_n  + u_np_i p_j+u_jp_ip_n)
+
z_3
\sinh\eta
 m_im_jm_n
$$

\s

\s

$$
+
L
 (m_jp_ip_n+m_n p_ip_j+m_i p_jp_n)
+
 r_2
\sinh\eta
\sin\theta
p_ip_jp_n
\Biggr].
$$
We obtain

{

\s

$$
H^2F^2C^i{}_{pq}C_{ijn}
=
-(
  u_6u_pu_q
+
L_3
m_p m_q
+
L_2
p_p p_q
)
(
  u_6u_ju_n
+
L_3
m_j m_n
+
L_2
p_j p_n
)
$$

\s

\s

$$
-
\Bigl(
L_2
( u_q p_p+u_pp_q)
+
L
 (m_pp_q+m_q p_p)
+
 r_2
\sinh\eta
\sin\theta
p_pp_q
\Bigr)
    \times
 $$

\s

\s

$$
\Bigl(
L_2
(u_n p_j+u_jp_n)
+
L
 (m_jp_n+m_n p_j)
+
 r_2
\sinh\eta
\sin\theta
p_jp_n
\Bigr)
$$

\s

\s

$$
-
\Bigl(
L_3
 (u_qm_p  +  u_pm_q)
+
z_3
\sinh\eta
m_pm_q
+
L
  p_pp_q
\Bigl)
    \times
 $$

\s

\s

$$
\Bigl(
L_3
 (u_nm_j  +  u_jm_n)
+
z_3
\sinh\eta
m_jm_n
+
L
 p_jp_n
\Bigr),
$$
after which we get

{

\s \s

$$
H^2F^2(C^i{}_{pq}C_{ijn}-C^i{}_{pn}C_{ijq})
$$

\s

\s

$$
=
-(
  u_6u_pu_q
+
L_3
m_p m_q
+
L_2
p_p p_q
)
(
  u_6u_ju_n
+
L_3
m_j m_n
+
L_2
p_j p_n
)
-[qn]
$$

\s

\s

$$
-
L_2
 u_q p_p
\Bigl(
L_2
(u_n p_j+u_jp_n)
+
L
 (m_jp_n+m_n p_j)
+
 r_2
\sinh\eta
\sin\theta
p_jp_n
\Bigr)
-[qn]
$$

\s

\s

$$
-
L_2
 u_p p_q
\Bigl(
L_2
(u_n p_j+u_jp_n)
+
L
 (m_jp_n+m_n p_j)
+
 r_2
\sinh\eta
\sin\theta
p_jp_n
\Bigr)
-[qn]
$$

\s

\s

        $$
-
L
m_pp_q
\Bigl(
L_2
(u_n p_j+u_jp_n)
+
L
 (m_jp_n+m_n p_j)
+
 r_2
\sinh\eta
\sin\theta
p_jp_n
\Bigr)
-[qn]
$$

\s

\s

        $$
-
L
m_q p_p
\Bigl(
L_2
(u_n p_j+u_jp_n)
+
L
 (m_jp_n+m_n p_j)
+
 r_2
\sinh\eta
\sin\theta
p_jp_n
\Bigr)
-[qn]
$$

\s

\s

        $$
-
 r_2
\sinh\eta
\sin\theta
p_pp_q
\Bigl(
L_2
(u_n p_j+u_jp_n)
+
L
 (m_jp_n+m_n p_j)
\Bigr)
-[qn]
$$

{

\s

\s

$$
-
L_3
u_qm_p
\Bigl(
L_3
 (u_nm_j  +  u_jm_n)
+
z_3
\sinh\eta
m_jm_n
+
L
 p_jp_n
\Bigr)
-[qn]
$$

\s

\s

\s

$$
-
L_3
u_pm_q
\Bigl(
L_3
 (u_nm_j  +  u_jm_n)
+
z_3
\sinh\eta
m_jm_n
+
L
 p_jp_n
\Bigr)
-[qn]
$$

\s

\s

\s

$$
-
z_3
\sinh\eta
m_pm_q
\Bigl(
L_3
 (u_nm_j  +  u_jm_n)
+
z_3
\sinh\eta
m_jm_n
+
L
 p_jp_n
\Bigr)
-[qn]
$$

\s

\s

\s

$$
-
L
  p_pp_q
\Bigl(
L_3
 (u_nm_j  +  u_jm_n)
+
z_3
\sinh\eta
m_jm_n
\Bigr)
-[qn].
$$

{

Now we reduce similar terms, arriving at
the representation

\s

$$
H^2F^2(C^i{}_{pq}C_{ijn}-C^i{}_{pn}C_{ijq})
$$

\s

\s

$$
=
-(
  u_6u_pu_q
+
L_3
m_p m_q
+
L_2
p_p p_q
)
(
  u_6u_ju_n
+
L_3
m_j m_n
+
L_2
p_j p_n
)
-[qn]
$$

\s

\s

$$
-
(L_2)^2
(u_j u_q p_pp_n+ u_p u_n p_jp_q)
-
LL_2
 u_q p_p
 (m_jp_n+m_n p_j)
-[qn]
$$
\s

\s

\s

$$
-
LL_2
 u_p p_q
m_n p_j
-
L^2
(m_nm_pp_q p_j+m_qm_jp_n p_p)
-[qn]
$$

\s

\s

        $$
-
L
L_2
m_pp_q
u_n p_j
-[qn]
-
LL_2
m_q p_p
(u_n p_j+u_jp_n)
-[qn]
$$

\s

\s

$$
-
L_3
u_qm_p
\Bigl(
L_3
  u_jm_n
+
L
 p_jp_n
\Bigr)
-
L_3
u_pm_q
\Bigl(
L_3
u_nm_j
+
L
 p_jp_n
\Bigr)
-[qn]
$$

\s

\s

$$
-
Lz_3
\sin\eta
(m_pm_q p_jp_n+m_jm_n p_pp_q)
-
LL_3
  p_pp_q
 (u_nm_j  +  u_jm_n)
-[qn]
$$

{

\nin
which can be written in the simpler form:

\s

$$
H^2F^2(C^i{}_{pq}C_{ijn}-C^i{}_{pn}C_{ijq})
$$

\s

\s

$$
=
-(
  u_6u_pu_q
+
L_3
m_p m_q
+
L_2
p_p p_q
)
(
  u_6u_ju_n
+
L_3
m_j m_n
+
L_2
p_j p_n
)
-[qn]
$$

\s

\s

$$
-
(L_2)^2
(u_j u_q p_pp_n+ u_p u_n p_jp_q)
-[qn]
$$

\s

\s

$$
-
LL_2
\Bigl(
 u_q p_p   m_jp_n
+
 u_p p_q
m_n p_j
+
m_pp_q
u_n p_j
+
m_q p_p
u_jp_n
   \Bigr)
-[qn]
$$

\s

\s

        $$
-
L^2
(m_nm_pp_q p_j+m_qm_jp_n p_p)
-[qn]
$$

\s

\s

$$
-
(L_3)^2
(u_ju_qm_p  m_n+u_pu_nm_j  m_q)
-[qn]
$$

\s

\s

$$
-
Lz_3
\sin\eta
(m_pm_q p_jp_n+m_jm_n p_pp_q)
-[qn]
$$

\s

\s

\be
-
LL_3
\Bigl(
  p_pp_q  (u_nm_j  +  u_jm_n)
+
p_jp_n(u_qm_p +u_pm_q )
\Bigr)
-[qn].
\ee
From the obtained expansion (B.4) the validity of the equalities (B.1)-(B.3) just follows.

{

\s \s

\def\bibit[#1]#2\par{\rm\noindent\parskip1pt
                     \parbox[t]{.05\textwidth}{\mbox{}\hfill[#1]}\hfill
                     \parbox[t]{.925\textwidth}{\baselineskip11pt#2}\par}

\nin {  REFERENCES}

\s

\bibit[1]  Asanov,G.S. (2017).  Pseudo-Finsleroid Metrics With Two Axes, \it Europian Journal of Mathematics. \rm
 DOI 10.1007/s40879-017-0160-6.

\bibit[2]
Asanov,G.S. (2015). Pseudo-Finsleroid metric function of spatially anisotropic relativistic type.  arXiv:1512.02268 [math.GM].

\bibit[3]  Minkowski,H. (1909). \it Raum und Zeit,\rm  Phys. Z. {\bf 10}, 104.

\bibit[4] Minkowski,H. (1915). \it  Das Relativit{\"a}tsprinzip,  \rm Ann. Physik. {\bf 47} , 927.

\bibit[5] Synge, J.L. \it Relativity: The General Theory. \rm  North-Holland, I960, 505pp.

\bibit[6] Cartan,E. \it  Les espaces de Finsler.  \rm  Actualit\'es  79,  \rm Paris, 1934, 77pp.

\bibit[7]  Rund,H. \it The Differential Geometry of Finsler  Spaces. \rm Springer, 1959, 284 pp.




\bibit[8]  Asanov,G.S. \it Finsler Geometry, Relativity and Gauge  Theories. \rm D.~Reidel Publ. Comp., 1985, 370 pp.

\bibit[9] Bao,D.,      Chern,S.-S.,  Shen,Z. \it An Introduction in Riemann-Finsler Geometry.\rm
Graduate Texts in Mathematics, vol. 200, Springer, New York  (2000).

\bibit[10]
{\it A Sampler of Riemann-Finsler Geometry.}
Edited by   Bao, D.,   Bryant, R.L.,   Chern,S. ,  Shen,Z.
Cambridge University Press, 2004, 376 pp.

\bibit[11] Matveev,V.S., Rademacher,H.-B.,  Troyanov,M., Zeghib,A. (2008).
Finsler Conformal Lichnerowicz-Obata conjecture,
arXiv:0802.3309 [math.DG].

\bibit[12] Matveev,V.S., Rademacher,H.-B.,  Troyanov,M., Zeghib,A. (2009).
Finsler Conformal Lichnerowicz-Obata conjecture,
\it Annales de l'institut Fourier, \bf 59(3), \rm 937-949.

\bibit[13]  Matveev,V.S.,  Troyanov,M. (2011).
The Binet–Legendre Metric in Finsler Geometry.
arXiv:1104.1647 [math.DG].

\bibit[14]  Matveev,V.S.,  Troyanov,M. (2012).
The Binet–Legendre Metric in Finsler Geometry
\it Geometry and Topology \bf 16,  \rm  2135-2170.

\bibit[15] Vincze,Cs. (2015). \it Average methods and their applications in differential geometry. \rm
I.J.Geom.Phys. \bf 92, \rm 194-209.

\end{document}